\documentstyle[eqsection]{article}

\newcommand{\Om}{\Omega}
\newcommand{\la}{\langle}
\newcommand{\ra}{\rangle}

\newenvironment{pf}{\noindent{\sc Proof}.\enspace}{\rule{2mm}{2mm}\medskip}
\newenvironment{pfn}{\noindent{\sc Proof} \enspace}{\rule{2mm}{2mm}\medskip}
\newtheorem{theorem}{Theorem}[section]
\newtheorem{proposition}{Proposition}[section]
\newtheorem{lemma}{Lemma}[section]

\newtheorem{remark}{Remark}[section]
\newtheorem{remarks}{Remark}[section]
\newtheorem{definition}{Definition}[section]
\newcommand{\be}{\begin{equation}}
\newcommand{\ee}{\end{equation}}

\newcommand{\wk}{\rightharpoonup}

\newcommand{\om}{\omega}

\newcommand{\ov}{\overline}

\newcommand{\wtilde}{\widetilde}
\newcommand{\R}{\bf R}

\newcommand{\T}{\bf T}

\renewcommand{\d }{\delta }

\newcommand{\e }{\varepsilon }

\newcommand{\vphi}{\varphi }

\newcommand{\dps}{\displaystyle}

\begin{document}

\title{{\bf Periodic solutions of 
nonlinear wave equations with general nonlinearities}}

\author{Massimiliano Berti, Philippe Bolle}
\date{}
\maketitle

{\bf Abstract:}
We prove the existence of 
small amplitude periodic solutions, with strongly irrational 
frequency $ \om $ close to one, for completely resonant
nonlinear wave equations. 
We provide multiplicity results 
for both monotone and nonmonotone nonlinearities.
For $ \om $ close to one we prove the existence of a large
 number $ N_\om $ of 
$ 2 \pi \slash \om $-periodic in time solutions
$ u_1, \ldots, u_n, \ldots, u_N $:
$ N_\om \to + \infty $ as $ \om \to 1 $.
The minimal period of the $n$-th solution $u_n $ is proved to be
$2 \pi \slash n \om $.
The proofs are based on a Lyapunov-Schmidt reduction and 
variational arguments.
\footnote{Supported by M.U.R.S.T. Variational Methods and Nonlinear
Differential Equations.}
\\[2mm]
Keywords: Nonlinear wave equation, Infinite dimensional 
Hamiltonian systems, 
variational methods, periodic solutions, Weinstein Theorem.
\\[1mm]
2000AMS subject classification: 35L05, 37K50, 58E05.

\section{Introduction}

In this paper we look for periodic solutions
of the completely resonant nonlinear wave equation
\be\label{eq:main}
\cases{ u_{tt} - u_{xx} + f ( u ) =  0, \cr
u ( t, 0 )= u( t,  \pi )  = 0.}
\ee
The wave equation (\ref{eq:main}) is said completely resonant
if $ f(0) = f'(0) = 0 $. Precisely,
we assume that  the nonlinearity $ f $ belongs to
$ {\cal C}^r ({\bf R}, {\bf R}) $,
for some large $ r $, and  satisfies
\begin{itemize}
\item
$(F1)$ $\  f (0) = f'(0) = \ldots = f^{(p-1)} (0) = 0, $
$ \ f^{(p)} (0) := a p! \neq 0 \ $
for some $ p \in {\bf N} $, $ p \geq 2 $.
\end{itemize}

Equation (\ref{eq:main}) can be viewed as an
{\it infinite dimensional Hamiltonian system}
with Hamiltonian
\be\label{Hamilt}
{\cal H} ( u, p ) :=  \int_0^\pi dx \ \frac{p^2}{2}
+ \frac{u_x^2}{2} + F(u)
\ee
where $ F(u) $ is a primitive of $ f $.
\\[1mm]
\indent
For finite dimensional Hamiltonian systems
existence of periodic orbits
has been widely investigated starting from the
pioneering work of Poincar\'e \cite{Po}.
The first result on the existence of periodic orbits
close to an elliptic equilibrium
is the well-known Lyapunov Center Theorem, see \cite{Ly}.
Assume that a Hamiltonian vector field possesses at $0$
an elliptic equilibrium
and let $\{ \om_j \}_{j=1, \ldots, n} $
be the frequencies of the linear oscillations close to $ 0 $.
Lyapunov's Theorem states that,
if one frequency, say $\om_i $, is
non resonant with the other ones,
i.e. $ \om_i l - \om_j \neq 0 $,
$ \ \forall l \in {\bf N}, j \neq i $, then the
periodic solutions of frequency $ \om_i $ of the linearized system
at $ 0 $ can be continued
into small amplitude periodic
solutions of the non-linear system with frequencies
$ \om $ close to $ \om_i $.
This result is proved using bifurcation theory
(see for example \cite{AP}, pag 145 and follows)
which can be applied thanks to the previous
non-resonance condition, yielding the existence of
smooth families of periodic solutions.

The existence of periodic orbits close to an
elliptic equibrium whose linear frequencies
$\{ \om_j \}_{j=1, \ldots, n} $ are possibly resonant,
but only at an order greater than $ 4 $,
i.e. $ \sum_{j=1}^n \om_j k_j \neq 0 $, $\forall
0 < |k| = \sum_{j = 1}^n |k_j| \leq 4 $, has been later proved
in the famous Birkhoof-Lewis Theorem \cite{BL}, see also \cite{Mo1}.
The proof is based on a normal form argument, an Implicit Function Theorem
and variational arguments.

Existence of small amplitude periodic
solutions of Hamiltonian systems with
a finite number of degrees of freedom,
  in case of  general resonances between the
linear frequencies, has been proved much later in
the celebrated Weinstein Theorem \cite{We}
(see also \cite{Mo} and \cite{Bar}).
The Weinstein Theorem
ensures, that, if the Hamiltonian $ H $ possesses at $0$ a non-degenerate
minimum, then, on each small energy level, there exist
at least $ n $ periodic orbits geometrically distinct.
This Theorem is proved through a
Lyapunov-Schmidt reduction: one first partially solves the above problem
by the techniques of bifurcation theory and then restricts the action
functional to the pseudo-solutions so obtained. It turns out that
the restricted functional is defined on a
manifold  homeomorphic
to\footnote{Here $S^{2n-1}:= \{ x \in {\bf R}^{2n} \ | \ |x|=1 \}.$}
$S^{2n-1} / S^1 $,
and  admits, for topological reasons, at least $n$
distinct critical points giving rise to periodic solutions
of the system. In general
such periodic solutions will be no more organized into smooth families.
Existence of periodic orbits with fixed period
close to an elliptic resonant equilibrium
have been later proved by Fadell and Rabinowitz in \cite{FR}.
The result of \cite{FR} is still based on a Lyapunov-Schmidt reduction
together with
a cohomological index to estimate  the number of
critical points of a reduced action functional (and hence of small
amplitude periodic solutions with fixed period). No information
on the minimal periods of these solutions is given.
\\[1mm]
\indent
It is natural to try to extend such type of results to
Hamiltonian PDE's as the wave equation, the plate equation,
the nonlinear Schr\"odinger equation, etc... .
The main difficulty in carrying over this program is the fact that,
in infinite dimensional systems,
``small denominator problems'' appear, whereas in finite dimensional systems
``small denominators'' appear only in the construction of quasi-periodic
solutions. This is simply explained as follows. The
frequencies $ \{ \om_j  \}_{j \in {\bf N}} $ of the linear oscillations
close to the equilibrium
typically grow as $ \om_j \approx j^d $ for some $ d > 0 $
and then the sequence
$\{ \om_i l - \om_j \}$, $ \ l \in {\bf N}, j \in {\bf N}, j \neq i $,
accumulates to $ 0 $.

Despite this difficulty, Kuksin \cite{K2} and Wayne \cite{W1},
were able to find, extending in a suitable way
KAM techniques, periodic (and even quasi-periodic)
solutions in some Hamiltonian PDE's in one spatial dimension
under  Dirichlet boundary conditions. As usual in KAM-type results,
the periods of such persistent solutions satisfy a strong
irrationality condition, as the classical Diophantine condition,
so that these orbits
exist only on energy levels belonging to some Cantor set
of positive measure; see for example \cite{Ku}
for a complete and updated exposition of the KAM approach.
The main limitation of this method is the fact that
standard KAM-techniques require the linear frequencies $ \om_j  $
to be well separated
(non resonance between the linear frequencies).

To overcome such difficulty
a new method for proving the existence of small amplitude
periodic solutions, based on the Lyapunov-Schmidt reduction,
has been developed in \cite{CW}.
Rather than attempting to make a series of canonical transformations
which bring the Hamiltonian into some normal form,  the
solution is constructed directly. Making the ansatz that a periodic solution
exists one writes this solution as a Fourier series and
substitutes that series into the partial differential equation.
In this way one is reduced
to solve two equations: the so called $(P)$ equation, which is
infinite dimensional, where small denominators appear, and the finite
dimensional $(Q)$ equation,
which corresponds to resonances. Due to the presence
of small divisors
the $(P)$ equation is solved by a Nash-Moser Implicit Function Theorem.
Later on, this method has been improved
by Bourgain to show the persistence of periodic solutions
in higher spatial dimensions \cite{Bo1} and of quasi-periodic solutions
in Hamiltonian perturbations of 2D linear Schr\"odinger
equations \cite{B3}.

A simple proof of the Lyapunov center Theorem in PDE's,
still based on the Lyapunov Schmidt procedure,
has been given more recently in \cite{B1}, where, making use of
a non-resonance condition on the frequency
stronger than the classical Diophantine one, the $(P)$ equation
is solved
by means of the standard Contraction Mapping Theorem.
Apart its simplicity, \cite{B1} has also the advantage of
avoiding some further condition present in \cite{CW}.
Other related results on the existence of periodic
solutions in nonlinear plate equations in higher dimension
have been given in \cite{BP2}.

We remark that in all the previous cases the $(Q)$ equation
is still finite dimensional since one handles cases with only
a finite number of resonances among the linear frequencies.

The first results on the existence of small amplitude
periodic solutions for some
{\it completely resonant} PDE's as (\ref{eq:main}) 
have been given in \cite{LS}, for the
specific nonlinearity  $ f ( u ) = u^3 $, and in \cite{BP1} when
$f(u)=u^{3}+ {\rm h.o.t.} $.
One looks for small amplitude periodic solutions of  (\ref{eq:main})
bifurcating from the infinite dimensional linear space of
solutions
$$
V :=
\Big\{ v_{tt} - v_{xx} = 0 , \ v( t , 0) = v( t, \pi) = 0, \
v( t + 2 \pi , x ) = v( t, x ) \Big\}.
$$
Here the linear frequencies  are
$ \om_j = j \in {\bf N} $.
The approach of \cite{BP1}
is still based on the Lyapunov-Schmidt reduction.
The $(P)$ equation is solved, for
the strongly irrational frequencies $ \om \in W_\gamma $,  where
$$
W_\gamma := \Big\{ \om \in {\bf R} \ \Big| \
| \om l - j | \geq \frac{\gamma }{l}, \ \forall l \neq j \Big\},
$$
through the Contraction Mapping Theorem
(it is proved in \cite{BP1} that, for $ 0<\gamma < 1 \slash 3 $,
the set $ W_\gamma $ is uncountable
and accumulates to $ \om = 1 $ both from
the left and from the right).
Next, the $(Q)$ equation, infinite dimensional,
is solved by looking for
non degenerate critical points of a suitable
functional (see remark \ref{nondeg})
and continuing them, by means of the Implicit Function Theorem,
into families of periodic solutions of the nonlinear equation.
Clearly the main difficulty in applying this method
is to check that the fore mentioned
non-degeneracy condition holds. This is verified in \cite{BP1}
for the particular case
$ f (0) = f'(0) = f'' (0) = 0 $ and $ f^{'''}(0) \neq 0$ for which
informations of the minimal periods of the solutions
are also given.
\\[1mm]
\indent
The aim of the present paper is to prove, assuming
only that the nonlinearity $f$ satisfies $(F1)$,
the existence of a large number
of small amplitude periodic solutions of (\ref{eq:main})
with fixed period,
see  Theorems \ref{thm:main1}, \ref{thm:main2}, \ref{theven},
\ref{thm:main3}.
Roughly, we show that equation (\ref{eq:main}) has a large number $ N_\om $ of
$2 \pi \slash \om $ periodic
solutions $ u_1, \ldots, u_n, \ldots, u_N $, where
$ N_\om \approx \sqrt{\gamma^\tau \slash (\om -1) } \to + \infty $ as
$ \om \to  1 $, for some $\tau \in [1,2]$. We expect the number of
small amplitude periodic solutions of frequency
$ \om $ to be in general $ O( \sqrt{1 \slash (\om -1) } ) $,
see remark \ref{rem:opt}.
We also prove that the minimal period of the $n$-th periodic
solution $ u_n $ is $ 2 \pi \slash n \om $.

These Theorems can also be considered a generalization of the results
of \cite{FR} (which hold for finite dimensional Hamiltonian systems) 
for the nonlinear wave equation (\ref{eq:main}) satisfying $(F1)$.

Theorems \ref{thm:main1}, \ref{thm:main2}, \ref{theven} and
\ref{thm:main3} are proved
using the Lyapunov-Schmidt procedure which leads to solve the
equations (P) and (Q).
We overcome the small denominator problem arising in the (P) equation,
by imposing on the frequency $ \om $ the same
strong irrationality condition as in \cite{BP1}. On the contrary,
the infinite dimensional $(Q)$ equation is solved by means of
a variational principle which provides periodic solutions of
(\ref{eq:main}) with fixed frequency: this is
the main novelty of our approach with respect to \cite{BP1}.
This variational principle is inspired by the perturbation method
in critical point theory introduced in \cite{AB}-\cite{ACE}.
With respect to these last two papers, the new feature of our approach
is the presence of an infinite dimensional $(Q)$ equation.

Let us give a more precise outiline of our arguments.
Critical points of the action functional
$\Psi ( u ) :=$
$\int_0^{2\pi} dt \int_0^{\pi} dx  \
(\om^2 \slash 2)  u_t^2 - (1 \slash 2 ) u_x^2 - F(u)$
defined in a suitable space of functions
$2 \pi $ periodic in time and vanishing at $x = 0 $ and $x = \pi $
are (weak) solutions of our problem.
First, for any $ v \in V $, we find
as in \cite{BP1}, a unique $ w(v) \in V^\bot $,
$ w(v) = O(||v||^p) $,
which solves the (P) equation. Then, according to the
variational procedure of \cite{AB}-\cite{ACE},
see also \cite{FR},
all the critical points of the action functional
$ \Psi $ close to the origin can be obtained as critical points of
the reduced action functional $ \Phi (v) := \Psi (v + w(v)) $,
defined on the linear space $ V $ by
$$
\Phi (v) := \Phi_\e (v)
= \frac{\e}{2} || v ||^2  - \int_0^{2 \pi} \ dt \int_0^\pi \ dx \
F(v + w(v)) - \frac{1}{2} f(v + w(v)) w(v) \qquad {\rm where } \qquad
\e := \frac{ \om^2 - 1 }{2},
$$
see lemma \ref{natcons} and remark \ref{rem:vic}.
By $(F1)$, $ \Phi_\e $ possesses a local minimum at the origin.
To show that $ \Phi_\e $ satisfies the geometrical hypotheses of
the Mountain Pass Theorem \cite{AR} the nonlinearity $ f $ plays its
role.
Here we have to make a different study whether the first
term in the Taylor expansion at $ 0 $
of the nonlinearity $ f( u ) = a u^p + h.o.t. $ is an odd power of $ u $,
or an even power of $ u $.
In the first case one has that
$$
\Phi_\e (v)
= \frac{\e}{2} || v ||^2
  - \frac{a}{ p+1} \int_0^{2 \pi} \int_0^\pi v^{p+1} + \ {\rm h.o.t.}
  $$
and hence $ \Phi_\e $ (or $- \Phi_\e$) satisfies the mountain pass
geometry provided that sign $ \e $=sign $a$. In other words,
when  $ a > 0$ (resp. $a < 0 $) we will find solutions for $ \om > 1 $
(resp. $\om < 1 $)
and the bifurcation is supercritical (subcritical).

The second case ($p$ even) is more difficult since
$ \int_0^{2 \pi} \int_0^\pi v^{p+1} \equiv 0 $ and then one must
look more carefully in the expansion of $ \Phi_\e $.
We point out that this difficulty is similar to the one
arising in the
Yamabe problem, see in \cite{AM}-\cite{BM},
where the first order expansion term in the corresponding reduced functional
identically vanishes and one needs a sharper analysis
for finding its critical points.
Let us consider for simplicity $ f(u) = u^p $ with $p $ even.
Then it is possible to show that
$$
\Phi_\e (v) = \frac{\e }{ 2} ||v||^2
- \frac{a^2 }{2} \int_\Omega v^p L^{-1} v^p + \ {\rm h.o.t.}
$$
where $z := L^{-1} v^p $ is a solution of $ z_{tt} - z_{xx} = v^p $,
$ z (t,0) = z(t,\pi)=0 $, which is $2\pi$ periodic in time,
see formulas (\ref{Geven})-(\ref{DReven}).
Now our existence proof requires a careful study
of the nonquadratic term $ (- a^2 \slash 2)$ $
\int_\Omega v^p L^{-1} v^p $, see lemmas
\ref{lemexpl}, \ref{cor:nl}, \ref{Gpos}. In particular, we show in
lemma \ref{Gpos} that
$ (- a^2 \slash 2)
\int_\Omega v^p L^{-1} v^p \geq 0$, and this implies that
$- \Phi_\e $ satisfies  the Mountain pass geometry
for $\e < 0$, i.e. $\om < 1$.

Next, we prove that the Palais-Smale sequences at
the mountain pass critical level converge to a critical
point which lies inside the domain of
definition of the functional (where the Lyapunov-Schmidt reduction
is defined) provided $|\e |$ is small enough,
i.e. $ \om $ is sufficiently close to $ 1 $.
This is done by means of suitable estimates
on the mountain pass critical level.

Furthermore, in order to find multiple critical points of $ \Phi_\e $,
we do not use
Ljusternik-Schnirelmann theory as for the symmetric
Mountain pass Theorem \cite{AR}
(remark that $ \Phi_\e $ is even).
Instead, we look for
``local'' mountain-pass critical points
of $ \Phi_\e $ restricted to the
subspaces $ V_n \subset V $ consisting by the functions of $ V $
which are $ 2 \pi \slash n $-periodic
w.r.t. time. Since each $ V_n $ is invariant for the
gradient flow of $ \Phi_\e $, any critical point of
$ {\Phi_\e}_{|V_n} $ is a critical point of $ \Phi_\e $
on the whole $ V $.
This method also  provides precise
informations on the minimal periods, on the
norm and on the energies of the solutions, see also remark
\ref{rem:nen}.
It also follows that,
as the frequency $ \om $ tends to $ 1 $, the number of
small amplitude solutions we find increases to $ + \infty $.

Let us further point out that
the preceding arguments would also improve the
results of \cite{FR} for finite dimensional Hamiltonian systems, 
giving informations on the
minimal periods of the solutions.

Finally we prove that the weak solutions found are
indeed classical $C^2$ solutions of (\ref{eq:main}).
\\[1mm]
\indent
We underline that our results do not require any monotonicity
assumption on the nonlinearity.
The major difficulty outlined above for proving the existence
of periodic solutions when the first term in the Taylor expansion
of the nonlinearity $ f $ is an even power of $ u $ reflects
an intrinsic difficulty of the problem and it is
encountered by all the known methods when
dealing with non monotone nonlinearities.
Up to our knowledge
the only existence result (in any case quite different in nature from ours)
for periodic solutions
without the assumption of monotonicity of the nonlinearity is given in
\cite{Co} where the existence of one, not small, periodic solution
with rational frequency and symmetry properties,
is proved through global variational techniques.
\\[1mm]
\indent
Before concluding this introduction we would like to  mention that
the first breakthrough in the study of bifurcation
of periodic solutions for nonlinear wave equations
was the work of Rabinowitz \cite{R0} (see also \cite{DST})
where the existence of forced periodic solutions under a monotonicity
assumption on the nonlinearity is proved.
For a global bifurcation result see \cite{R1}.
Moreover, large free vibrations
have been found by
Rabinowitz in \cite{R} (see also \cite{BCN})
using global variational methods. Both in \cite{R}
and \cite{BCN} the nonlinearity is assumed to be
(superquardatic and) monotone. In these papers the period of the
solutions of (\ref{eq:main}) is a rational multiple of $\pi$.
Actually this condition introduces in the problem
strong compactness properties, while in the case of
irrational periods,
``small-denominator phenomenona'', similar to the ones
described above, would appear into the proof.
\\[2mm]
\indent
The paper is organized as follows. In section 2 we perform the
Lyapunov-Schmidt reduction and we define the variational formulation.
In section 3 we prove a general abstract result on existence of
critical points. In section 4 we apply the former result
to the nonlinear wave equations. In section 5 we prove the
regularity of the solutions.
\\[3mm]
{\bf Acknowledgments:}
The authors thank Prof. A. Ambrosetti, Prof. D. Bambusi and
Dr. M. Procesi for useful discussions.
Part of this paper was written when the second
author was  visiting S.I.S.S.A. in Trieste.

\section{The variational setting}

We look for periodic solutions of (\ref{eq:main}) 
with frequency $ \om $ 
as $ u ( t, x ) = q ( \om t , x) $, where $ q( \cdot , x)$  
is $ 2 \pi $ periodic. We 
are then led to solve the problem (still denoting $ q \equiv u $)
\be\label{eq:freq}
\om^2 u_{tt} - u_{xx} + f( u ) = 0, \qquad u(t,0)= u(t, \pi) = 0, \qquad
u( t+ 2 \pi, x )= u(t,x).
\ee
We will look for solutions of (\ref{eq:freq})
$ u : \Omega \to {\bf R} $  where $ \Omega 
:= {\bf R} \slash 2 \pi {\bf Z} \times (0, \pi )$,   
in the following Banach space 
$$
X := \Big\{  u \in H^1 ( \Omega, {\bf R} ) 
\cap L^\infty ( \Omega, {\bf R} )\ |  \ u(t,0)= u(t, \pi) = 0, 
\ u(-t,x) = u(t,x) \Big\} 
$$
endowed, for $ \om \in [ 1/2, 3/2 ] $ , $ \om \neq 1 $, with the norm
\be \label{defnorm}
| u |_{\om} := |u|_{\infty}+ |\om -1|^{1/2} ||u||_{H^1}.
\ee
Note that we can restrict to the space $ X $ of functions even in time
because equation (\ref{eq:main}) is reversible. 

Any $ u \in X $ can be developed in Fourier series as
$$
u(t,x) = \sum_{l \geq 0, j \geq 1} u_{lj} \cos lt \sin jx, 
$$
and its $H^1$ norm and scalar product are written as 
$$
||u||^2 :=||u||_{H^1}^2 :=  
\int_0^{2\pi} dt \int_0^\pi dx \ u_t^2 + u_x^2
= \frac{\pi^2}{2} \sum_{l \geq 0, j \geq 1} 
u_{lj}^2 ( j^2 + l^2 ),
$$
\be\label{scalpro}
( u, w ) := \int_\Omega \ dt \ dx \ 
u_t w_t + u_x w_x 
= \frac{\pi^2}{2} \sum_{l \geq 0, j \geq 1} u_{lj} w_{lj} 
( l^2 + j^2 ) \qquad  \forall u , w \in X.
\ee
We will also denote with 
$ |u|_{L^2} := ( \int_\Om |u|^2 )^{1/ 2 } =   
( \pi \slash \sqrt{2} ) ( \sum_{l \geq 0, j \geq 1} u_{lj}^2 )^{1/ 2 }$ 
the $L^2$-norm of $u\in $ and with $( u, w )_{L^2} := \int_\Omega u w $ the 
 $L^2$-scalar product.

We need the following  preliminary lemma proved in the Appendix.

\begin{lemma}\label{composition}
The Nemitski operator $ u \to f(u) $ 
is well defined and $ C^1 $ on $ X $. Its differential at $ u $ is given by
$ h \to Df(u)[h] = f'(u) h $. Moreover there is $ \rho_0>0 $ and
positive constants $ C_1, C_2 $, depending only on $ f $,
such that, $ \forall u \in X $ with
$ |u|_\om < \rho_0 $, 
\be\label{nemitski}
| f(u) |_{\om} \leq  C_1 |u|_\infty^{p-1} |u|_\om \leq C_1 |u|_\om^p,
\ee
\be\label{nemitski2}
|f'(u)h|_\om \leq C_2 |u|_\infty^{p-2} |u|_\om |h|_\om 
\leq C_2 |u|_\om^{p-1} |h|_\om.
\ee
\end{lemma}

The  solutions in $ H^1_0( \Omega , {\bf R} ) $ 
of the linear equation $ v_{tt} - v_{xx} = 0 $  that are even in time
are the elements of
$$
V := \Big\{ v(t,x) = \sum_{j \geq 1} \xi_j \cos ( j t ) \sin j x \ 
\Big| \ \xi_j \in {\bf R}, \  \sum_{j\geq 1} j^2 \xi_j^2 <
+\infty  \Big\}.
$$
The space $ V $ can be also written as
$$
V := \Big\{ v(t,x) = \eta ( t + x ) - \eta ( t - x ) \ 
\Big| \  \eta ( \cdot ) \ \in H^1 ({\bf T}), \ 
\eta {\rm \ odd} \Big\}
$$
where ${\bf T} := {\bf R} \slash 2 \pi {\bf Z}$ is the 
one dimensional torus and 
$ \eta(s) = \sum_{j \geq 1 } ( \xi_j / 2 ) \sin (j s) $. 
Note that 
$ \la \eta \ra := (1 / 2 \pi ) \int_{{\bf T}} \eta (s) \ ds = 0 $.

$ V $, endowed with the
$H^1$-norm $|| \ ||$, is a Hilbert space, 
continuously embedded in $ X $, since, 
for $ v =\sum \xi_j \cos(jt) \sin(jx) \in V$, we have 
$ |v|_\infty \leq \sum_{j\geq1} |\xi_j| \leq C ||v||$
by the Cauchy-Schwarz inequality. Moreover the embedding
$ (V , || \cdot ||) \to (V, | \cdot |_\infty ) $ is compact. 
\\[1mm]
\indent
In order to get multiplicity results we will consider in the sequel 
the subspaces of $V$
\be\label{V_n}
V_n := \Big\{ v \in V \ \Big| \ v \ \ { \rm is } \ \ \frac{2\pi}{n}-
{\rm periodic} \ \ {\rm w.r.t.} \ \ t \Big\},
\ee
for $ n \geq 1 $. 
For $ v = \eta( x+t ) - \eta( x-t ) \in V $ and any $ n \in {\bf N} $ 
we define 
$$
({ \cal L}_n v) (t,x) := \eta (n(t+x))-  \eta (n(t-x)). 
$$  
It is immediately realized that 
$$
V_n =  { \cal L}_n V =  
\Big\{ v (t,x) = \eta (n(x+t)) - \eta (n(x-t))  \ \Big| \   
\eta ( \cdot ) \in H^1 ({\bf T}), \ \eta ( \cdot ) {\rm \ odd } \Big\}.
$$

We next consider 
the Lagrangian action functional $ \Psi : X \to {\bf R}$
defined by 
\be\label{acfunct}
\Psi(u) :=
\int_0^{2\pi} dt \int_0^{\pi} dx  \  
\frac{\om^2}{2} u_t^2 - \frac{1}{2} u_x^2 - F(u)
\ee
where $ F ( \cdot ) $ is the primitive of  $f$ such that $F(0)=0$, i.e.
$ F(u) = \int_0^u f(s) \ ds$. It is easy to see that 
$ \Psi \in C^1 ( X, {\bf R}) $ and that 
\be
D \Psi (u) [h] = \int_0^{2\pi} dt \int_0^{\pi} dx  \  
\om^2 u_t h_t  - u_x  h_x  - f(u) h, \quad \forall h \in X.    
\ee
Critical points of $ \Psi $ are weak solutions of (\ref{eq:freq}). 
\\[1mm]
\indent
In order to find critical points of $ \Psi $  
we perform a Lyapunov-Schmidt reduction.
The space $ X $ can be decomposed as  
$ X = V \oplus W $ where
$$
W := 
\Big\{ w \in X \ \Big| \ (w,v)_{L^2} = 0, \ \forall \ v\in V 
\Big\}=
\Big\{ \sum_{l \geq 0, j \geq 1} w_{lj} \cos ( l t ) \sin j
x \in X  \ 
\Big| \ w_{jj} = 0 \ \forall  j \geq 1 \Big\}.
$$ 
In fact, if $ u=\sum u_{lj} \cos (lt) \sin (jx) \in X$, then the 
$H^1$ norm of $\Pi_V (u):=\sum u_{jj} \cos (jt) \sin(jx)$ is
finite, so that $\Pi_V (u) \in X$. As a consequence $ \Pi_W (u):=
u - \Pi_V (u) \in X $. Moreover the projectors $ \Pi_V : X \to V $ and
$ \Pi_W : X \to W $ are continuous. 
$W$ is also the $ H^1 $-orthogonal of $V$ in $X$. Note that 
if $ v = \Pi_V u $ has minimal period w.r.t time $ 2 \pi $ 
then also $ u $ has minimal period w.r.t. time $ 2 \pi $. 
\\[1mm]
\indent
Setting $ u := v + w $ with $ v \in V $ and $ w \in W $, 
(\ref{eq:freq}) is equivalent to the following system of 
two equations (called resp. the $(Q)$ and the $(P)$ equations) 
$$
(Q) \qquad
- \om^2 v_{tt} + v_{xx} = \Pi_V f(v + w), 
$$
$$
(P) \qquad -\om^2 w_{tt} + w_{xx} = \Pi_W f(v + w). 
$$
Note that $w$ solves the $(P)$ equation (in a weak sense) iff it is a
critical point of the restricted functional 
$ w \to \Psi (v + w) \in {\bf R}$, {\it i.e.} if and only if  
\be\label{critW}
D \Psi(v + w) [h] =   
\int_0^{2\pi} dt \int_0^{\pi} dx  \ \om^2 
w_t h_t  - w_x h_x - f(v+w) h = 0, \quad \forall  h \in W .
\ee
We will first solve the $(P)$ equation in lemma \ref{lem:red}  
through the standard Contraction Mapping Principle,
assuming that $ \om \in { \cal W } := \cup_{ \gamma > 0 } W_\gamma $ 
where $ W_\gamma $ is the  set of strongly non-resonant frequencies
introduced in \cite{BP1}
\be
W_\gamma := \Big\{ \om \in {\bf R} \  \Big| \
|\om l - j | \geq \frac{\gamma}{ l }, \quad \forall j \neq l \Big\}. 
\ee
As proved in remark [2.4] of \cite{BP1}, for $ \gamma < 1 \slash 3 $,
the set $ W_\gamma $ is uncountable,
has zero measure and accumulates
to $ \om = 1 $ both from the left and from the right. 
It is also easy to show that $ W_\gamma = \emptyset $
for $ \gamma \geq 1 $.

The proofs of the next lemmas \ref{lem:operL} and \ref{lem:red} are given
in the Appendix.

\begin{lemma}\label{lem:operL} 
For $ \om \in W_\gamma $, let us define 
$$
\e =\frac{\om^2-1}{2} \quad (\e \sim \om-1 \ \ {\rm as} \ \ \om
\to 1).$$  
Then   the operator 
$ L_\om := -\om^2 \partial_{tt} + \partial_{xx} : D(L_\om) \subset
W  \to W$ has a
bounded inverse $ L_\om^{-1} : W \to W $ which writes
\be\label{Lominv}
L_\om^{-1} h := \sum_{l \geq 0, j \geq 1, j \neq l} 
\frac{h_{lj}}{ (\om^2 l^2 - j^2) } \cos (lt) \sin (jx)
\ee
and satisfies, for a positive constant $ C_3 $ independent of
$\gamma$ and $ \om $, 
 $ | L_\om^{-1}\Pi_W u|_\om  \leq (C_3 \slash \gamma) |u|_\om$ for
$ u \in X $.  
Let $ L^{-1} : W \to W $ be the inverse operator of 
$- \partial_{tt} + \partial_{xx} $. 
There exists $C_4 > 0$ such that $\forall r,s \in X$ 
\be \label{diffLLom}
\Big|\int_{\Omega} r(t,x) \Big( L_{\om}^{-1}-L^{-1} \Big)
(\Pi_W s(t,x)) \ dt \ dx \Big| \leq C_4 \frac{\e}{ \gamma } 
|r|_{\om} |s|_{\om}.
\ee
Moreover, $\forall r,s \in X$
\be \label{Lomest}
\Big|\int_{\Omega} r(t,x)  L_{\om}^{-1}
(\Pi_W s(t,x)) \ dt \ dx \Big| \leq C_4 \Big( 1 + \frac{\e}{ \gamma } \Big) 
|r|_{\om} |s|_{\om},
\ee
and finally, $\forall s \in X$,
\be \label{estsupp}
|L_\om^{-1} \Pi_W s|_{L^2} \leq C_4 
\Big( 1 + \frac{\sqrt{\e}}{ \gamma } \Big) |s|_{\om}.
\ee
\end{lemma}
 \begin{lemma} \label{lem:red}
There exists $ \rho > 0 $ such that,  $ \forall v \in 
{\cal D}_\rho := \{ v \in V \ | \ |v|_\om^{p-1} \slash \gamma < \rho \}$
there exists a unique $ w(v) \in W $ with $|w(v)|_\om <  |v|_\om$ 
solving the $(P)$ equation. Moreover, for some positive
constant $ C_5 $, we have 
\begin{itemize}
\item
$ i) \ |w(v)|_\om \leq C_5 |v|_\om^p / \gamma$; 
$\ |w(v)|_{L^2} \leq C_5 (1+ \sqrt{\e}/ \gamma) |v|_\om^p $; 
\item
$ ii) \ $ for all $r \in X$, 
$$ \Big| \int_\Om \Big( w(v)  - L_\om^{-1} \Pi_W f ( v ) \Big) \  r \Big| \leq
\frac{C_5}{\gamma} \Big( 1+ \frac{\e}{\gamma} \Big) |v|_\om^{2p-1} |r|_\om; $$
\item
$ iii) \   w(-v)(t,x) = w(v)(t + \pi, \pi -x ) $; 
\item
$ iv)$ $ v \in V_n \ \Rightarrow \ w(v)$ is $2\pi/n$ periodic
w.r.t. $t$;
\item
$ v) $ \  the map $ v \to w(v)$ is in  $ C^1(V, W) $. 
\end{itemize}
\end{lemma}

In order to solve the $(Q)$ equation 
we consider  
the reduced Lagrangian action functional 
$ \Phi_\e : {\cal D}_\rho \to {\bf R} $ defined by
\be\label{redfun1}
\Phi_\e ( v ) := \Psi( v + w ( v ) ) =
\int_0^{2\pi} dt \int_0^{\pi} dx \ 
\frac{\om^2}{2} (v_t + (w(v))_t)^2 - \frac{1}{2} (v_x + (w(v))_x)^2 -
F(v + w(v)).
\ee
By (\ref{redfun1}) and formula (\ref{critW}) with $ h = w(v) $ 
we can write the reduced action functional as
\be\label{redfun2}
\Phi_\e ( v ) = \int_{\Omega} \  dt \ dx \  
\frac{\om^2}{2} v_t^2 - \frac{v_x^2}{2} - F(v + w(v)) +
\frac{1}{2} f(v + w(v)) w(v),
\ee
Since $ |v_t|_{L^2}^2 = | v_x |_{L^2}^2 = || v ||_{H^1}^2 / 2 $, we have 
\be\label{expfunc}
\Phi_\e (v)    
= \frac{\e}{2} || v ||^2  + \int_{\Omega} \  dt \ dx \
 \frac{1}{2} f(v + w(v)) w(v)-F(v + w(v)).
\ee
Defining the linear operator ${\cal I} : X \to X $
by $({\cal I} u) (t,x) = u(t+\pi,\pi-x)$ it results that
$\Psi \circ {\cal I}= \Psi $. Hence,
by lemma \ref{lem:red}-$iii$) 
and since $ -v = {\cal I} v $,  
functional $ \Phi_\e $ is even:  
\be\label{eveness}
\Phi_\e ( - v ) = \Psi (-v+w(-v))=\Psi ({\cal I}
(v+w(v))) = \Psi(v+w(v))=\Phi_\e ( v ), \qquad \forall 
v \in {\cal D}_\rho,
\ee
By lemma \ref{lem:red}-$v)$ functional 
$ \Phi_\e $ is in $ C^1 (V, {\bf R}) $.  
Arguing as in \cite{AB}-\cite{ACE}-\cite{FR} 
we find that, $ \forall h \in V $ 
\be\label{diffred}
D \Phi_\e (v) [ h ] = D \Psi (v+w(v))[h] = 
\int_\Omega dt \ dx \ \om^2 v_t h_t  -  v_x h_x  - f(v + w(v))h
= \e (v,h) - \int_\Omega dt \ dx \ f(v + w(v))h.
\ee
Hence the following lemma holds 

\begin{lemma}\label{natcons}
If $ v \in V $ is a critical point of the reduced action functional 
$ \Phi_\e : {\cal D}_\rho \to {\bf R} $ then $ u = v + w(v)$ is a
critical point of $ \Phi_\e : X \to {\bf R}$. 
\end{lemma}

\begin{remark}\label{rem:vic}
On the other hand it is possible to prove that any critical point 
$ u $ of $ \Phi_\e $, sufficiently close to $ 0 $, 
can be written as
$ u = v + w $ with $ w=w(v) \in W, v \in V $ and $v $ is a 
critical point of $ \Phi_\e $. 
\end{remark}

\section{An abstract result}

Let $ E $ be a Hilbert space with scalar product $( \cdot, \cdot )$ 
and norm $|| \cdot ||$. We will denote $B := \{ v \in E \ | \ ||v|| \leq 1 \}$
the closed unit ball of $ E $
and  $ S := \{ v \in E \ | \ ||v||=1 \} $ the unit sphere.
\\[1mm]
\indent
Let us consider a  $ C^1 $ functional
$\Phi : B_{r_0} \subset E \to {\bf R}$, 
defined on the ball $ B_{r_0} := \{ v \in E \ | \ ||v|| \ < {r_0} \} $, by 
\be\label{absfunc}
\Phi (v)= \frac{\mu}{2} ||v||^2 - G(v) + R (v)
\ee
where $ \mu \in (0, + \infty ) $ and 
\begin{itemize}
\item 
{\bf $(H1)$} \  
$ G \in C^1( E, {\bf R}) $ is positively homogeneous of degree $ q + 1 $
with $ q > 1 $,
i.e. $ G(\lambda v) = \lambda^{q+1} G(v) $ $ \forall \lambda \in {\bf R}_+ $,  
and $ DG : V \to V^* $ is compact;
\item {\bf $(H2)$} \ 
$ R \in C^1( B_{r_0}, {\bf R}) $ and $ DR : B_{r_0} \to V^* $ is compact;
\item {\bf $(H3)$} \  $ R(0) = 0 $ and 
$ \exists \alpha > 0 $ such that 
$ | DR (v) [v] | \leq \alpha  \ ||v||^{q+1},$ 
$ \forall v \in B_{r_0}.$
\end{itemize}
\noindent
By the compactness of $ DG : V \to  V^* $, $ G $ maps 
the unit sphere $ S $
into a bounded subset of $ {\bf R} $ (indeed 
$ G(v) = \int_0^1 DG(sv)[v] ds $ and $DG(B)$ is bounded in $ V^* $).
It follows that 
the functional $ U : E \backslash \{ 0 \} \to { \bf R } $
defined by 
$$
U(v) := \frac{G(v)}{||v||^{q+1}},
$$
which satisfies
$ U ( \lambda v ) = U ( v ) $, $\forall  v \in E \backslash \{ 0 \}$, 
$\forall \lambda \in {\bf R_+}$,  is bounded.
We shall assume the further condition 
\be\label{positive}
m := \sup_{ v \in E \backslash \{ 0 \} } U(v) > 0 .
\ee
The following lemma, quite standard, 
ensures that $ U $ attains its maximum on $ S $.

\begin{lemma}
$K_0:= \{ v  \in S \ | \  G (v) = m \} = \{ v  \in S \ | \  U (v) = m \}$ 
is not empty and compact.
\end{lemma}

\begin{pf} Since $DG(B)$ is compact,
$\forall \d >0$ there is a finite dimensional subspace $F_\d$ of 
$E$  such that
$$
\forall v \in B \  \ \forall h \in F_\d^\bot \ \
|DG(v)[h]| \leq \d ||h||.
$$ 
As a result, calling $ P_\d $ the orthogonal projection onto $ F_\d $, 
we obtain for $v\in B$
$$
|G(v)-G(P_\d v)| 
\leq \int_0^1 |DG((1-s)P_\d v +s v)[v-P_\d v]| \ ds
\leq \int_0^1 \d ||v-P_\d v|| \ ds  \leq \d ||v-P_\d v|| \leq \d.
$$ 
One can derive from this property that $G_{|B}$ is continuous
for the weak topology in $E$. Indeed $G_{|B} $ is the uniform limit,
as $\d \to 0 $, of  $ G \circ P_\d $ and each function 
$ G \circ P_\d $
is continuous in the weak topology since $ P_\d $ 
is a linear compact operator. Since $B$, the unit
ball of $ E $, is compact for the weak topology, $G$
attains its maximum on $B$ at some point $v_0 \in B$.
We have $G ( v_0 ) = m > 0 $, so that $ v_0 \neq 0 $, and
$ U ( v_0 ) = m / ||v_0||^{q+1} \geq m$. Hence, by the definition
of $m$,  $U(v_0)=m$ and $v_0 \in S$, {\em i.e.} $v_0\in K_0$.

To prove that $K_0$ is compact, consider a sequence $(v_k)$ with 
$v_k \in K_0 $ for all $ k $. Since $(v_k)$ is bounded, we may
assume that (up to a subsequence) $(v_k)$ converges to some $ v \in E $
for the weak topology. For all $ k $,  $ v_k \in B$  hence,
since $B$ is closed for the weak topology, 
$v \in B $.
Moreover, since $G_{|B}$ is continuous for the weak topology,
$G(v)=m$. By the same argument as above, we get that $|| v || = 1 $. 
Hence $(v_k) \rightharpoonup v $ and $(||v_k||) \to ||v||$ and we
can conclude that $(v_k) \to v$ strongly.   
\end{pf}

The main result of this section 
-which will allow us in the next section to get our 
multiplicity results- is the following Proposition.

\begin{proposition} \label{mainlem}
Assume that the functional $ \Phi \in C^1( E, {\bf R}) $ 
satisfies $(H1)$, $(H2)$, $(H3)$ and (\ref{positive}).  
There exists a small positive constant $ C_0 $, depending only 
on $ q $, such that, if
\be\label{estionep}
\frac{\alpha}{m} 
\leq C_0
\quad {\rm and} \quad
\Big( \frac{\mu}{m}  \Big)^{\frac{1}{q-1}} \leq C_0 {r_0},
\ee
functional $ \Phi $ has a critical point $ v \in B_{r_0} $
on a critical level  
\be \label{cepbound}
c = \frac{q-1}{2} m \Big( \frac{\mu}{(q+1)m} \Big)^{\frac{q+1}{q-1}} 
\Big[ 1 +  O \Big( \frac{\alpha}{m} \Big) \Big].
\ee
Moreover 
\be\label{estidist}
v = \Big( \frac{\mu}{m(q+1)} \Big)^{\frac{1}{q-1}} y \qquad 
{\rm with} \qquad
{\rm dist}( y , K_0) \leq  h \Big(\dps 
\frac{\alpha}{m}  \Big)
\ee
for some function $ h $, depending only on $ G $,
which satisfies $\lim_{s \to 0} h (s ) = 0 $.
\end{proposition}

\begin{remark}\label{nondeg}
The key difference with the approach of \cite{BP1} is the following:
instead of trying to show that the functional
$(\mu \slash 2 ) ||v||^2 - G(v) $ possesses non-degenerate 
critical points (if ever true), and then continuing them through the
Implicit Function Theorem 
as critical points of $\Phi $, 
we find critical points of $\Phi $  
showing that the nonlinear perturbation term $ R $   
does not affect the mountain pass geometry of the functional  
$ (\mu  \slash 2 ) ||v||^2 - G(v) $. 
Actually, in \cite{BP1}, non degenerate 
critical points of $ G $ constrained on $S$
are continued. This is equivalent to what said before since
any critical point $ \wtilde{v} $ of $ G $  
constrained to $ S $ gives rise, by homogeneity, to a 
critical point 
$(\mu \slash (q+1) G(\wtilde{v}))^{1\slash q-1} \wtilde{v} \in V $
of the functional $ (\mu  \slash 2 ) ||v||^2 - G(v) $.
\end{remark}

The remainder of this section is devoted to the proof of
Proposition \ref{mainlem}. Let us define the positive constant
$ \delta := [2 (q+1)]^{1/q-1}$. Assuming that 
\be\label{epsmall}
\Big( \frac{2\mu}{m} \Big)^{\frac{1}{q-1}} < {r_0}
\ee
it results that
the rescaled functional $ \varphi $   
\be\label{rescafun}
\varphi (y) := 
\frac{1}{(q+1)m} \Big( \frac{(q+1)m}{\mu} \Big)^{\frac{q+1}{q-1}}
\Phi \Big( \Big( \frac{\mu}{(q+1)m} \Big)^{1/q-1} y \Big)
\ee
is well defined on an open neighborhood of $ \ov{B}_\delta $.
We shall also assume that 
\be\label{etasmall}
\eta(\alpha,m):= \frac{\alpha}{(q+1)m}  \leq \eta_1.
\ee
(\ref{epsmall}) and (\ref{etasmall}) stem from  
(\ref{estionep}) if $C_0$ is chosen appropriately.
We shall prove that, for $ \eta_1 $ small enough, 
$ \varphi $ has a critical point $ y $
on a critical level $ d $ such that 
\be \label{cepbound2}
d = \Big( \frac{1}{2}-\frac{1}{q+1} \Big) (1+ O(\eta(\alpha,m)))
\ee
and 
\be\label{distbound}
{\rm dist} ( y , K_0 ) \leq 
h \Big( \frac{\alpha}{m} \Big) 
\ee
for some function $h$ such that $\lim_{s\to 0} h(s)=0$. 
(\ref{cepbound2}) and (\ref{distbound}) imply resp. 
(\ref{cepbound}) and (\ref{estidist}) since, by the definition of 
$\varphi$, 
$ v := \displaystyle \Big( \frac{\mu}{(q+1)m} \Big)^{1/q-1} y$
will then be a  
critical point of $ \Phi $  at  level 
$c = d (q+1)m $ $ \displaystyle 
\Big( \frac{\mu}{(q+1)m} \Big)^{q+1\slash q-1}. 
$ 
We can write
\be\label{S1}
\varphi (y) = \frac{||y||^2}{2} - \frac{G(y)}{(q+1)m} + S (y),
\ee
where, by the properties of $ R $,
\be \label{S2}
\forall y \in B_\delta, \ \ 
|DS (y)[y]| \leq \eta( \alpha, m ) ||y||^{q+1}
\leq  \eta_1 ||y||^{q+1}.
\ee
As a consequence, since $S(y)= \int_0^1 ds DS (sy)[y] $, 
\be\label{S3}
|S(y)|\leq  \frac{\eta(\alpha,m)}{q+1} \  ||y||^{q+1} 
\leq \frac{\eta_1}{q+1}  \ ||y||^{q+1}.
\ee

We need the following technical lemma:

\begin{lemma}\label{elem}
Let $W$ be some (small) neighborhood of $0$ in ${\bf R}$.
There are constants $  C_6 , C_7 >0$ such that,
$ \forall a \in [1/2,1]$, $\forall \  0 < \xi \leq C_6, \forall \ 
0< \eta \leq C_6 $, for
all $ C^1 $ function $ r : [0,\delta] \to {\bf R}$
with $ \sup_{t\in [0,\delta]} |r(t)|+|r'(t)| \leq $ $ \xi $, 
the map$f_{a,r} : [0,\d] \to \R$ defined by 
$ f_{a,r}(t) := (t^2/2)-a(t^{q+1}/q+1)+r(t)$, satisfies
\be \label{tt}
t\in [0,\delta] \backslash W \ , \ |f'_{a,r}(t)| \leq \eta 
\Longrightarrow |t- a^{-1/q-1}| \leq C_7 (\xi +\eta),
\ee
and
\be \label{lowbound}
\inf_{t\in [0,\delta]\backslash W, |f_{a,r}'(t)| \leq \eta} f_{a,r}(t) \geq 
\Big( \frac{1}{2}-\frac{1}{q+1} \Big) 
\Big( \frac{1}{a} \Big)^{2/q-1}-C_7(\xi+\eta).
\ee
\end{lemma}
\begin{pf}
Consider first the function $f_a (t) :=\frac{t^2}{2}-a
\frac{t^{q+1}}{q+1}$. $ f'_a $ possesses outside $W$ a unique zero  
$ t_a := a^{-\frac{1}{q-1}}$ at the level
$$
f_a(t_a)= \Big( \frac{1}{2}-\frac{1}{q+1} \Big) 
\Big(\frac{1}{a}\Big)^{\frac{2}{q-1}}.
$$
$ f''_a(t_a)=1-q < 0 $, and there are $ \nu >0$, $\ov{\nu}>0$ such that,
$ \forall a \in [ 1 / 2, 1 ]$,
\be \label{star}
\forall t \in [t_a-\nu, t_a+\nu] \  \
f''_a(t) \leq - \nu/2 \quad {\rm and} \quad 
\forall t \in [0,\delta] \backslash (W\cup [t_a-\nu, t_a+\nu])
\  \  |f'_a(t)| \geq \ov{\nu}.
\ee
Choose $C_6  > 0 $ such that $2 C_6 <\ov{\nu} $. If 
$ \sup_{t\in [0,\delta]}|r(t)|+|r'(t)| 
\leq \xi \leq C_6$ and $|f'_{a,r}(t)|
\leq \eta \leq C_6$ for some $ t \notin W$,  
then $|f'_a(t)|\leq \xi+\eta < \ov{\nu}$ and therefore,
by the second inequality in (\ref{star}), 
$ t \in [t_a-\nu, t_a+\nu]$. Moreover by $(\ref{star})$ again,
$| f_a'(t)| \geq (\nu \slash 2 )|t - t_a | $ for $t \in [t_a-\nu, t_a+\nu]$, 
and hence $ |t-t_a| \leq 2 (\xi +\eta)/
\nu $, which proves estimate $(\ref{tt})$. 
Finally, by the Taylor formula $ f_a(t)= f_a (t_a ) + O( ( \xi+\eta )^2)
$, so that $ f_{a,r}(t)=f_a(t_a)+ O(\xi)+ O((\xi+\eta)^2)$ and
estimate (\ref{lowbound}) follows straightforwardly.  
\end{pf}

Let us define the open set
$$
A := \Big\{ v \in  E \backslash \{0 \} \ \ \Big| \  
U(v) > \frac{m}{2} \Big\},
$$
the min-max class of path
$$
\Gamma = \Big\{ \gamma \in C ( [0,1],A \cup \{0\} ) \ \ \Big| \ \
\gamma ( 0 ) = 0, \ ||\gamma(s)|| \leq \delta, \  \ 
|| \gamma (1)|| = \delta \Big\}
$$
and the ``restricted-to-$A$'' mountain pass
level of $ \varphi $
$$
d := \inf_{\gamma \in \Gamma} \max_{s \in [0,1]}
\varphi (\gamma(s)).
$$
Before proving that $ d $ is a critical level for $\varphi$,
we shall check that, provided $\eta_1$ has been chosen small enough,
it satisfies $(\ref{cepbound2})$. 
For $ y \in S_A := S \cap A$ we define $f_y : [0,\delta]
\to {\bf R}$ by $f_y(t)=\varphi (ty)$. 
By (\ref{S1}) we have 
\be\label{eq:fy}
f_y(t)=\frac{t^2}{2}-\frac{U(y)t^{q+1}}{(q+1)m} +
r_y(t),
\ee
where, by (\ref{S3}) and (\ref{S2}),
\be\label{eq:ry}
|r_y(t)|=|S(ty)| \leq \frac{\eta(\alpha,m)}{q+1} t^{q+1}
\qquad {\rm and}  \qquad
|r'_y(t)| \leq \eta(\alpha,m) t^{q}. 
\ee
We first prove an upper bound for $d$.
For $ y_0 $ belonging to $ K_0 $ we clearly have
$ d \leq \max_{t\in [0,\delta]} f_{y_0}(t)$ and hence,
from (\ref{eq:fy}) and (\ref{eq:ry}), 
$$
d \leq \max_{t\in [0,\delta]} \Big( \frac{t^2}{2} -
\frac{t^{q+1}}{q+1} \Big) +
\frac{\eta(\alpha,m)}{q+1}\delta^{q+1}
= \Big( \frac{1}{2}-\frac{1}{q+1} \Big) +
\frac{\eta(\alpha,m)}{q+1}\delta^{q+1}.
$$

For the lower bound, we observe that there is a
small neighborhood $W$ of $0$ in ${\bf R}$ such that
$\forall y \in S_A $ and $\forall t\in W $, $f'_y (t) > 0$, whereas,
by the choice of $ \delta =[2 (q+1)]^{1 \slash q-1} $ 
and provided $\eta_1$ is small enough, 
$f'_y (\delta) < 0$. By the continuity 
of $D\varphi$, this implies that 
$\forall \gamma \in \Gamma $ there is $ s \in [0,1]$ 
such that $ D \varphi (\gamma (s))[ \gamma (s)] = 0 $,
$\gamma(s) \neq 0 $. As a consequence,
$$
d \geq \inf_{y \in S_A, f'_y(t)=0, t \in [0,\delta]\backslash  W} f_y(t)
\geq \inf_{y \in S_A, f'_y(t)=0, t \in [0,\delta]\backslash  W} 
\frac{t^2}{2} - \frac{t^{q+1}}{(q+1)} + r_y(t)
$$
and the desired lower bound follows  by  (\ref{eq:ry}) and lemma \ref{elem}.

In order to complete the proof of Proposition
\ref{mainlem}, we need the following

\begin{lemma} \label{mpass}
There exists a sequence $(y_k)$,
$y_k \in A\cap B_{\delta}$ such that either $\varphi (y_k) \to
d$ and $ D \varphi (y_k) \to 0 $, or
$\varphi (y_k) \to d$ and $ U(y_k) \to m / 2 $ and
  $D\varphi (y_k) [y_k] \to 0$. 
\end{lemma}

\begin{pf}
Arguing by contradiction, assume that no
sequence $(y_k)$ satisfies the above properties. 
Then there exists $\nu >0$ such that, for all $v\in A\cap B_\delta$, 
\be \label{ps1}
|\varphi (v) - d| < \nu \quad \Longrightarrow \quad 
||D\varphi (v)|| > \nu
\ee
and
\be \label{ps2}
|\varphi (v)- d|+ \Big| U(v)- \frac{m}{2} \Big| 
< \nu \quad \Longrightarrow \quad
|D\varphi (v)[v]| >\nu.
\ee
There exists a locally Lipschitz-continuous map 
$ X : A_\nu:=\{ v \in A \cap B_\delta \ | \ ||D\varphi (v)||>\nu \} \to E$ 
such that $ ||X(v)|| \leq 2||D\varphi (v)||$ and
$ D\varphi (v)[X(v)] \geq ||D\varphi (v)||^2 \geq 
\nu^2$ for all $v\in A_\nu$ ($X$ is called a pseudo gradient 
field for $\varphi$; for a proof of  existence, see for
example \cite{Str}).

Let $J : {\bf R} \to [0,1]$ be a smooth function such that
$J(s)=0$ if $|s| \geq \nu/2$ and $J(s)=1$ if $ |s| \leq \nu/4$.

Let $ \beta_1 $ and $ \beta_2 < \d $ be such that 
$ \forall v \in A $ with $ ||v|| \leq \beta_1 $ or $ || v || \geq \beta_2$,
$ \varphi (v) < d -\nu $. 
$ \beta_2 $ does exists because $ \forall y \in S_A$,
 by (\ref{eq:fy})-(\ref{eq:ry}),
$ \vphi( \delta y) \leq $ $( \d^2 \slash 2 ) - (\d^{q+1} \slash 2(q+1)) 
+ $ $ O(\eta (\alpha, m)\d^{q+1}  )$ $ =  \d^2 $ $ (-( 1 \slash 2)+ $ 
$ O (\eta (\alpha, m) )  ) \leq 0 $, by the choice of $\d$ 
provided $ \eta_1 $ is small enough.
Let $ K : {\bf R} \to [0,1] $ be smooth and satisfy $K(s)
=0$ if $|s| \leq \beta_1^2/4$ or $ |s| \geq (\beta_2+\delta)^2/4$
and $K(s)=1$ if $\beta_1^2 \leq |s| \leq \beta_2^2 $. We define
$$
Y(v)=K(||v||^2)  J(\varphi (v)- d) \Big[ - \Big( 
1-J \Big( U(v)- \frac{m}{2} \Big) \Big)  X(v) -
J \Big( U(v)- \frac{m}{2} \Big) \    (D\varphi (v)[v]) \  v \Big].
$$
The vectorfield $Y$ is well defined on $A' := A \cup B_{\beta_1}$ 
and vanishes everywhere in $B_{\beta_1/2}$ and outside 
$B_\d$. It is locally Lipschitz continuous. Moreover it is
bounded, and it is easy to see that its flow $\phi_r (v)$ is 
defined, for all
time $r$, in $A'$ (notice that in a neighbourhood of $\partial A$, $Y(v)
\in {\bf R} v$, so that the flow leaves $U(v)$ invariant, and
one does not reach $ \partial A $);  $r \mapsto \phi_r (v) $ is the 
unique solution of $(d \slash dr) \phi_r (v) = Y (\phi_r (v) ) $
with $ \phi_0 (v) = v $.
We observe also that if 
$ \gamma \in \Gamma $ then $ \gamma_r $ defined by
$\gamma_r (s)=\phi_r(\gamma(s))$ belongs to $\Gamma$
for all $ r $ because $\phi_r (y)=y $ for 
$ || y || \leq \beta_1/2 $ and $||y|| \geq \d $.

Note that if $J(\varphi(v)-d) \neq 0 $ then $|\varphi (v)-d|
< \nu /2$. In that case, either $|U(v)-m/2|<\nu /2$ and then
by (\ref{ps2}) $|D\varphi (v)[v]| \geq \nu$, or 
$|U(v)-m/2|\geq \nu /2$ and then $Y(v) = - K(||v||^2)  J(\varphi (v)- d)
X(v)$. In both cases, by (\ref{ps1}) and the properties of
$X(v)$, 
\be\label{-gradflow}
D\varphi (v) [Y(v)] \leq -K(||v||^2)  J(\varphi (v)- d) \nu^2
\leq 0.
\ee
Now, if $|\varphi (v)-d| < \nu /4$ then $J(\varphi (v)-d)=1$
and by the definition of $\beta_1,\beta_2$, $\beta_1 \leq ||v|| 
\leq \beta_2$, so that $K(||v||^2)=1$.  
Therefore, by (\ref{-gradflow}), for all $ v \in A' \cap \ov{B}_\delta$,
\be \label{estflow}
|\varphi (\phi_r(v)) - d| < \frac{\nu}{4} \quad \Longrightarrow \quad 
\frac{d}{dr} \varphi (\phi_r (v)) \leq -\nu^2
\ee
Pick some $ \gamma \in \Gamma $ which satisfies $\max_{s\in [0,1]}
\varphi (\gamma (s)) \leq d +\nu/4$. By 
$(\ref{estflow})$, 
$\max_{s\in [0,1]}
\varphi (\phi_{1/2\nu}(\gamma (s))) \leq d -\nu/4$,
which contradicts the definition of $ d $.
\end{pf}

We can now prove that $ d $ is a critical level.
Let $ y_k $ be the sequence given by lemma 
\ref{mpass}.   
If $ D \varphi (y_k) [y_k] \to 0 $ and $ U ( y_k ) \to  m / 2 $ 
then, writing $ y_k = t_k z_k $ with $ z_k \in S $, we have   
$ \lim_k f_{z_k}'(t_k) = 0 $ and $U(z_k) \to m \slash 2 $
($|| y_k || \geq \rho_0 > 0 $ since $ \vphi ( y_k ) \to d > 0 $).
Therefore in this case, by (\ref{eq:fy}), (\ref{eq:ry}) 
and lemma \ref{elem}, 
if $ \eta_1 $ is small enough, 
$$ 
\lim \inf_k \vphi(y_k) \geq 
\Big( \frac{1}{2}-\frac{1}{q+1} \Big) 
2^{2/q-1} + O( \eta_1 ) > d
$$
by  estimate (\ref{cepbound2}). It 
follows that $ D \varphi (y_k) \to 0$ and
it is standard to prove (by the compactness
properties of $ DG $ and $DS$) that, up
to a subsequence, $(y_k)$ converges strongly
to a critical point $ y $ of $ \varphi $ 
on the level $ d $ (Palais Smale property, see e.g. \cite{Str}).
Note that, since $\vphi (y_k) to d $ and $ y_k \in A $, 
$y_k \in B_{\beta_2} $ for all $k $ so that $ y \in B_\d $. 

Let us now prove the last claim of Proposition \ref{mainlem}.
Writing $y=t z$ with $z\in S$ and $f'_z(t)=0$, by lemma \ref{elem} we get
$$
d \geq \Big( \frac{1}{2} - \frac{1}{q+1} \Big) 
\Big( \frac{m}{G(z)} \Big)^{2/q-1} - C_7 \eta(\alpha,m)
(\delta^q + \delta^{q+1}/(q+1)),
$$
and hence, by the upper estimate on $d$,
$G(z)=m+O(\eta(\alpha,m))$. Moreover by
$(\ref{tt})$ we have
$$
t=\Big( \frac{m}{G(z)} \Big)^{1/q-1} 
+ O(\eta(\alpha,m))=1+O(\eta(\alpha,m)).
$$ 
Hence
${\rm dist }(y,K_0)\leq g(C_8 \eta(\alpha,m)) + C_8 \eta(\alpha,m)$
for some positive constant $C_8$, where 
$$
g(s) := \sup_{v\in S, G(v) \geq m-s} {\rm dist}(v,K_0), \  \ \ 
\lim_{s\to 0^+} g(s)=0.
$$
The last claim of Proposition \ref{mainlem} follows 
with $ h(s) = g(C_8 s) + C_8 s $. 

\section{Applications to nonlinear wave equations} \label{secappl}

In this section we prove  
existence and multiplicity results
of periodic solutions for the nonlinear wave equation (\ref{eq:main})
satisfying $(F1)$. 
Precise informations on their minimal periods will be also given.

We start with some preliminary lemmas.

\begin{lemma}
Let $ m : {\bf R}^2 \to {\bf R}$ be locally integrable and 
$2\pi$-periodic w.r.t. both variables. Then 
$$
\int_0^{2\pi} \int_0^{\pi} m(t+x,t-x) \ dt \ dx =
\frac{1}{2} \int_0^{2\pi} \int_0^{2\pi} m(s_1,s_2) \
ds_1 \ ds_2.
$$ 
\end{lemma}

\begin{pf} 
Using new variables $(s_2,x)=(t-x, x)$ and the periodicity
of $m$, we get
$$
\int_0^{2\pi} \int_0^{\pi} m(t+x,t-x) \ dt \ dx =
\int_0^{2\pi}   \int_0^{\pi} m(s_2+2x,s_2) \ dx \ ds_2 =
\int_0^{2\pi} \  ds_2 \int_0^{\pi} m(s_2+2x, s_2) \ dx.
$$
By the periodicity of $ m $, setting $ s_1 := s_2 + 2 x $, we get
$$
\int_0^{\pi} m(s_2+2x, s_2) \ dx = \frac{1}{2} 
\int_0^{2\pi} m(s_1, s_2) \ ds_1, 
$$
which yields the result.
\end{pf}

\begin{lemma}\label{lem:ort}
If $ v \in V$ then $ v^{2p} \in W$.
\end{lemma}
\begin{pf}
Writing $ v(t,x) = \eta (t+x) - \eta (t-x) \in V $ we get
$ \forall u ( t, x ) = q ( t + x ) - q ( t - x ) \in V $,
$$
\int_0^{2\pi} \int_0^{\pi} v^{2p} (t,x) u(t,x) \ dx \ dt =
\frac{1}{2} \int_0^{2\pi} \int_0^{2\pi} 
\Big( \eta(s_1)-\eta(s_2) \Big)^{2p} (q(s_1)-q(s_2)) \ ds_1 \ ds_2 =0,
$$ 
because $(s_1,s_2) \mapsto (\eta(s_1)-\eta(s_2))^{2p} (q(s_1)-q(s_2))$
is an odd function. This implies that $ v^{2p} \in W $.
\end{pf}

We now apply the arguments of Proposition \ref{mainlem}
to functionals $\Phi $ of the form (\ref{absfunc})
defined on the Hilbert space $ V $ endowed with
the $H^1$ norm.

\begin{lemma}\label{miniper}
Assume that $ \Phi \in C^1 ( B_{r_0}, {\bf R}) $ satisfies the 
hypotheses of Proposition \ref{mainlem} and 
furthermore that $ \forall n \geq 2$ and 
$ \forall v \in V $ for which $ G({\cal L}_n v) > 0 $,
\be\label{smallGn}
G ({\cal L}_n v ) < n^{q+1} G(v).
\ee
Then, there is a constant $ C_9 \leq C_0 $ depending only on $ G $ such that,
provided $\alpha /m \leq C_9$ and $(\mu/m)^{1/(q-1)} \leq C_0 r_0 $, the
critical point $v$ of Proposition \ref{mainlem} 
has minimal period (w.r.t. $t$) $ 2 \pi $.
\end{lemma}

\begin{pf}
By Proposition \ref{mainlem}
we know that $v = ( \mu / m(p+1))^{1/p-1} y $ with dist$(y,K_0) \leq
h(\alpha /m)$ and $ \lim_{s \to 0} h(s) = 0 $ ($h$ depending on $G$ 
only). 
The lemma is proved once we show that $ \forall n \geq 2 $, 
$ V_n \cap K_0 =  \emptyset $, i.e. any $ z \in K_0 $ has
minimal period $ 2 \pi $. Indeed, if this is true, since
$ K_0 $ is compact, there exists $ \xi > 0 $ 
such that any $ y \in V $ with
dist$(y,K_0) \leq \xi $ has minimal period $ 2\pi $. Hence  for 
$\alpha /m$ small enough so that $h(\alpha /m) < \xi $, 
the critical point $ v $ has minimal period $ 2 \pi $.

For $ z (t,x) = \eta(t+x)-\eta(t-x) \in V $ we have 
\begin{eqnarray*}
|| z ||^2 & = & \int_0^{2\pi} \int_0^{\pi} (\eta'(t+x)-\eta'(t-x))^2 \  
dt \ dx
= \frac{1}{2} \int_0^{2\pi}  \int_0^{2\pi} (\eta'(s_1)-\eta'(s_2))^2
\ ds_1 \ ds_2 \\
&=& \frac{1}{2} \Big[ 2\pi  \int_0^{2\pi} \eta'(s_1)^2 \ ds_1 - 2
\int_0^{2\pi} \eta'(s_1) \ ds_1  \int_0^{2\pi} \eta'(s_2) \ ds_2
+ 2\pi \int_0^{2\pi} \eta'(s_2)^2 \ ds_2 \Big] \\
&=& 2\pi \int_0^{2\pi} \eta'(s)^2 \ ds.
\end{eqnarray*}
Hence
$ ||{\cal L}_n z||^2 = $ $ 2\pi \int_0^{2\pi} n^2 \eta'(ns)^2 \ ds 
= n^2 ||z||^2, $ because $ \eta $ is $ 2 \pi $-periodic.
As a consequence (for $z \neq 0$),
$$
U({\cal L}_n z)= \frac{G({\cal L}_n z)}{||{\cal L}_n z||^{q+1}}=
\frac{G({\cal L}_n z)}{n^{p+1}||z||^{q+1}} < U(z),
$$ 
by (\ref{smallGn}).
Hence for all $v\in V_n \backslash \{ 0 \}$, $n\geq 2$,
$U(v)  < \max_{ z \in V \backslash \{ 0 \} } U ( z ) = m $ and therefore
$ V_n \cap K_0 = \emptyset$.
\end{pf}

The next lemma provides a suitable approximation for 
the functional $ \Phi_\e $ defined in (\ref{expfunc}).

\begin{lemma} \label{appro}
The functional $ \Phi_\e : {\cal D}_\rho \to {\bf R} $ defined in 
(\ref{expfunc}) can be developed as
$$
\Phi_\e(v)= \frac{\e}{2} ||v||^2 - \int_{\Omega} F(v)
-\frac{1}{2} \int_\Om f(v) L^{-1}\Pi_W f(v) + R_{\om,1}(v),
$$
with
$$
DR_{\om,1}(v)[v]=O \Big( \frac{\e}{\gamma} |v|_\om^{2p} +
\frac{ |v|_\om^{3p-1}}{\gamma}  \Big).
$$
\end{lemma}
\begin{pf} We have, since 
$ \int_\Omega f(v) L^{-1} \Pi_W f'(v)v = 
\int_\Omega f'(v)v L^{-1} \Pi_W f(v) $,  
\begin{eqnarray*}
DR_{\om,1}(v)[v]&=&\int_\Om \Big[ 
f(v)+f'(v)L^{-1} \Pi_W f(v) -f(v+w(v)) \Big] v
\\
&=& \int_\Om \Big[ f(v)+f'(v)w(v)-f(v+w(v)) \Big] v + \int_\Om 
\Big[ ( L^{-1}-L_\om^{-1} ) \Pi_W f(v) \Big] f'(v)v \\ 
& + & \int_\Om \Big[ L_\om^{-1}\Pi_W f(v)-w(v) \Big] f'(v)v \\
& = & O \Big( |v|_\om^{p-1} |w(v)|_{L^2}^2 \Big) + 
O \Big( \frac{\e}{\gamma}  |v|_\om^{p} |v|_\om^p\Big) + 
O \Big( \frac{(\e+\gamma)|v|_\om^{2p-1}}{\gamma^2} |v|^p_\om 
\Big) \\
&=&  O \Big( \frac{\e}{\gamma^2} |v|_\om^{3p-1} + 
|v|_\om^{3p-1} \Big) + 
O \Big( \frac{\e}{\gamma} |v|_\om^{2p} \Big) + 
O \Big( \frac{\gamma +\e}{\gamma^2} |v|_\om^{3p-1}  \Big)
=  O \Big( \frac{\e}{\gamma} |v|_\om^{2p} + 
\frac{|v|_\om^{3p-1}}{\gamma} \Big)
\end{eqnarray*}
by lemmas \ref{composition},  \ref{lem:red}, formula (\ref{diffLLom})
and because $|v|_\om^{p-1} / \gamma \leq \rho$.
\end{pf}
 
Critical points of the restriction of $ \Phi_\e $ 
to the subspace $ V_n \subset V $  are  
critical points of $ \Phi_\e $. 

\begin{lemma}\label{lem:critn}
A critical point $ v_n $ of $ \Phi_{\e|V_n} :  V_n \cap {\cal D}_\rho
\to { \bf R } $ 
is a critical point of $ \Phi_\e : V \to {\bf R} $. 
\end{lemma}

\begin{pf}
We know that $ D \Phi_\e ( v_n ) [h] = 0 $, 
$ \forall h \in V_n $ and we want to prove that 
$ D \Phi_{\e} ( v_n ) [h] = 0 $, $ \forall h \in V_n^\bot \cap V $.
By (\ref{diffred}) it is sufficient to prove that 
$ \forall h \in V_n^\bot \cap V $
$$
\int_0^{\pi} dx \ \int_0^{2 \pi} \ dt \
f(v(t,x) + w(v)(t,x)) h(t,x) = 0.
$$
This holds true because,
by lemma \ref{lem:red}-$iv$),
$ f( v(t, x) + w(v)(t,x)) $ is $ 2 \pi \slash n $ 
periodic w.r.t $ t $
and $ t \to h(t,x) \in V_n^\bot \cap V $ does not contain any harmonic
with a frequency multiple of $ n $.
\end{pf}

Let us define $ \Phi_{ \e,n} : \{ v \in V \ | \ {\cal L}_n v 
\in {\cal D}_\rho \} \to {\bf R} $ by
$$
\Phi_{\e,n} (v) := \Phi_\e ({\cal L}_n v).
$$
By lemma \ref{lem:critn}, if $v$ is a critical point
of $ \Phi_{\e,n} $ (of minimal period $2\pi$), then $v_n := {\cal L}_n v$ 
is a critical point of $\Phi_\e$ (of minimal period $2\pi/n$). Hence, by lemma 
\ref{natcons}, $ u = v_n + w(v_n) $ is a solution of (\ref{eq:freq})  
and it has minimal period with respect to time $2\pi/n$.

In the sequel we shall use the norm 
\be\label{normn}
|v|_{\om,n} := |{\cal L}_n v|_\om =|v|_\infty + \e^{1/2} n ||v||.
\ee
Note that if $ | \e n^2 | < 1 $, which 
we shall always assume in the sequel, then 
$|v|_{\om,n} \leq C ||v|| $. 
As a consequence there is $ \rho_0 > 0 $ such that
$ {\cal D} := \{ v \in V \ | \ ||v||^{p-1} \slash \gamma \leq \rho_0 \}
$ is included in the domain of definition of $ \Phi_{\e,n} $.

\bigskip

Now we must distinguish between different cases depending whether the first
term in the Taylor expansion of 
of the nonlinearity $ f( u ) $ at $ 0 $ is an odd power of $u$
(like e.g. $ f ( u ) = a u^3 + h.o.t. $) or an even power of $u$
(like e.g. $ f ( u ) = a u^2 + h.o.t. $). 
We consider first the easiest
case of an odd dominant power for $ f(u) $. 

\subsection{Case I: $p$ odd}

In this case, by lemma \ref{appro}, 
the reduced action functional $ \Phi_\e $ defined in (\ref{expfunc})
can be written in the form (\ref{absfunc}) with
\be\label{Gpodd}
\mu=\e=\frac{\om^2-1}{2}, \quad  \quad
G(v) := a \int_{\Omega} \frac{v^{p+1}}{p+1}, \qquad \qquad 
a = \frac{f^{(p)} (0)}{p!} \neq 0
\ee
($G$ is homogeneous of degree $ p + 1 $) and  
$$
R(v) := R_{\om,1} ( v ) 
- \frac{1}{2} \int_\Om f(v) L^{-1}\Pi_W f(v) - 
\int_\Om \Big( F(v) - a \frac{v^{p+1}}{p+1} \Big)
$$
satisfies
\be\label{restopodd}
|DR (v)[v]| \leq C \Big( |v|_\om^{p+3} + \frac{\e}{\gamma}
|v|_\om^{2p} + \frac{|v|_\om^{3p-1}}{\gamma} \Big).
\ee 
In fact $ \int f(v)v - a v^{p+1} = O(|v|^{p+3}_\om) $
since,  by lemma \ref{lem:ort}, 
$ \int v^{p+2} = 0 $. Moreover $\int f'(v)v L^{-1} \Pi_W f(v) =
O(|f'(v)v|_{L^2}|f(v)|_{L^2})= 
O(|v|_\om^{2p})=O(|v|_\om^{p+3})$ since $p\geq 3$.
\\[1mm]
\indent
If $ a > 0 $, the functional 
$ G ( \cdot ) $ defined in (\ref{Gpodd}) is positive 
($p+1$ is even) and then
condition (\ref{positive}) is verified. 
On the contrary, if $ a < 0 $, 
we will look for critical points of $ - \Phi_\e $ and we will
take $ \e < 0 $, i.e. $ \om < 1 $. 
For definiteness we will assume
in the sequel that $ a > 0 $ and so $ \e >0 $.
\\[1mm]
\indent
We can easily prove that 
\begin{lemma}\label{poddl}
$ DG :V \to V^* $ and $ DR  : {\cal D} \to V^* $ 
are compact maps. 
\end{lemma}

\begin{pf}
Let $ || v_k || $ be bounded.
Then, up to a subsequence, $ v_k \wk \ov{v} \in V $
weakly in $ H^1 $ and $ v_k \to \ov{v} $ in $ | \cdot |_\infty $. 
Moreover since $ w_k := w(v_k) $ too is bounded we can also assume 
that $ w_k \wk \ov{w} $ weakly in $ H^1 $ and
$ w_k \to \ov{w} $ in $ | \cdot |_{L^q} $ norm for all
$ q < \infty $.  
Since $ DG(v)[h] = \int_\Omega a v^p h $, 
it follows easily that $ DG( v_k ) \to DG(\ov{v}) $ in $V^*$. 
We claim that $ DR  ( v_k ) \to \ov{R}  $ where 
$ \ov{R}[h] = \int_\Omega f(\ov{v} + \ov{w})h - a \ov{v}^p h
 =  DR ( \ov{v} ) [h ]$. 
Indeed, since $ w_k \to \ov{w} $ in $| \cdot  |_{L^q} $, it converges
(up to a subsequence) also $ a.e $. We can deduce, by the 
Lebesgue dominated convergence theorem, that
$ f ( v_k + w_k) \to f( \ov{v} + \ov{w} )$ in $ L^2 $
since  $f ( v_k + w_k) \to f( \ov{v} + \ov{w} )$ a.e.
and $ (f ( v_k + w_k)) $ is bounded in $ L^{\infty} $.
Hence, since 
$ D R (v_k) [h] = \int_\Omega f(v_k + w_k)h - a v_k^p h $,
$ D R (v_k) \to \ov{R}$.
\end{pf}

$ G(v) $ is nonnegative for all $ v \in V $ 
and it is easy to see that 
\be\label{Gnodd}
G( {\cal L}_n (v) ) = G(v), \quad \forall v \in V,  \ n \in {\bf N} . 
\ee
Hence 
$$
\Phi_{\e,n}(v) :=
\Phi_{\e} ( {\cal L}_n v) = \frac{\e n^2}{2} ||v||^2 - G(v) + R_n(v), \qquad 
\forall v \in {\cal D}
$$
where $ R_n(v) := R({\cal L}_n v)$. Of course $ DR_n $ (as $ DR $) is
compact. Moreover, by (\ref{restopodd}) and (\ref{normn})
\be \label{estDRn}
|DR_n (v)[v]|=|DR ({\cal L}_n v)[{\cal L}_n v]| \leq C \Big( |v|_{\om,n}^{p+3} + \frac{\e}{\gamma}
|v|_{\om,n}^{2p} + \frac{|v|_{\om,n}^{3p-1}}{\gamma} \Big) 
\leq C \Big( ||v||^{p+3} + \frac{\e}{\gamma}
||v||^{2p} + \frac{||v||^{3p-1}}{\gamma} \Big).
\ee

We want  to  apply Proposition \ref{mainlem}
in order to get  the existence, for $\om $ sufficiently close to $ 1 $,
of a critical point of $\Phi_{\e,n}$. 
With the notations of Proposition \ref{mainlem}, we set
$$
{r_0}=\frac{1}{C_0} \Big( \frac{\e n^2}{m}  \Big)^{1/p-1}.
$$
It is clear that $ B_{r_0} \subset {\cal D} $ provided 
$\e n^2 / \gamma $ is small.

Moreover, if $ \e n^2 \slash \gamma $ is small enough,  
the hypotheses of Proposition 
\ref{mainlem} are satisfied by $\Phi_{\e,n}$. In fact the first condition 
in (\ref{estionep}) is satisfied,   
because, by (\ref{estDRn}), 
$$
\alpha=O \Big( {r_0}^2 + \frac{\e}{\gamma} {r_0}^{p-1} +
\frac{{r_0}^{2p-2}}{\gamma} \Big) = O \Big( (\e n^2)^{2/(p-1)} + \Big(
\frac{\e n^2}{\gamma} \Big)^2 \Big)
$$
(here $ q = p $ and $\mu = \e n^2 $). The second condition 
in (\ref{estionep}) is trivially satisfied by the choice of $ r_0 $. 

Let us define, for $ \omega \in  { \cal W } := \cup_{\gamma>0}W_\gamma $, 
$ \gamma_\om := \max \{ \gamma \ | \ \om \in W_\gamma \}$.
By (\ref{Gnodd}), condition (\ref{smallGn}) is satisfied. 
As a conclusion, 
using Proposition \ref{mainlem} and lemma \ref{miniper} we obtain the
following result (we obtain pairs of solutions because $ \Phi_\e $ 
is even; if $u(t,x) $ is a solution of (\ref{eq:main}) then 
so is  $ u( t + \pi, \pi - x) $ ) 

\begin{theorem}\label{thm:main1}
Let $f $ satisfy $(F1)$ for an odd integer $ p \geq 3 $.  
Then there exists a positive constant 
$ C_{10} := C_{10} ( f ) $ such that, $ \forall \om \in {\cal W}$ and
$ \forall n\in {\bf N}\backslash \{0\} $   satisfying  
\be\label{estionep1}
\frac{| \om -1 | n^2}{\gamma_\om} \leq C_{10}
\ee
and $ \om > 1 $ if $ a > 0 $ (resp. $ \om <1 $ if $ a  < 0 $),
equation (\ref{eq:main}) possesses at least one pair of even
periodic (in time) classical $ C^2 $ 
solutions of minimal period $2\pi / (n \om )$.
\end{theorem}

The regularity of the solutions $u_n $ is proved in section \ref{sec:reg}.
By condition (\ref{estionep1}) we can find 
$ N_\om \approx O(\gamma_\om \slash \sqrt{|\om - 1|})$ periodic solutions
of (\ref{eq:main}) with period $ 2 \pi \slash \om $.
For a discussion on the optimality of the estimate (\ref{estionep1})
see remark \ref{rem:opt}. 

\begin{remark}
The unique small amplitude 
$ 2 \pi $ periodic in time solution of (\ref{eq:main}) is the equilibrium
solution $ u \equiv  0 $. Indeed it must be a critical point of the 
functional $ \Phi_0 (v) = (a/ p+1) \int v^{p+1} + R(v) $. 
Using (\ref{restopodd}) we see that its unique critical point is $ v = 0 $.
\end{remark}

\subsection{Case II: $p$ even}

The case in which $ f $ satisfies $(F1)$ for some even
$ p $ is more difficult to deal with since,
by lemma \ref{lem:ort}, $\int_\Omega v^{p+1} = 0 $, $ \forall v \in V$.
Then one must look more carefully for the dominant nonlinear term
in the reduced functional $ \Phi_\e $. 
We must distinguish different cases. Let $ f $ satisfy 
$ ( F1 ) $ for some even $ p $ and    
\begin{itemize}
\item {\bf $(N1)$} 
$ f^{(d)} ( 0 ) = b d ! \neq 0 $ for some odd integer 
$ d < 2 p - 1 $ and 
$ f^{(r)} ( 0 ) = 0 $ for any odd integer $ r < d $; 
\item {\bf $(N2)$} 
$ f^{(r)} ( 0 ) = 0 $ for any odd integer $ r  \leq 2 p -1 $;
\item {\bf $(N3)$} 
$f^{(r)} ( 0 ) = 0$ for any odd integer $ r < 2 p - 1 $ and  
$ f^{(2p-1)} ( 0 ) =: b (2p-1)! \neq 0$. 
\end{itemize}

If $(N1)$ holds
the dominant term in the reduced functional 
is supplied by the odd nonlinearity $b u^d$ in the
Taylor expansion of $ f $ at $0$  and then
one reduces immediately to the situation
discussed in the previous subsection. The following Theorem holds

\begin{theorem}\label{thm:main2}
Let $f $ satisfy $(F1)$ with $ p $ even and (N1). 
Then there exists a positive constant 
$ C_{11} := C_{11} ( f ) $ such that,  $ \forall \om \in {\cal W}$ and
$ \forall n \in {\bf N}\backslash \{0\} $ satisfying  
\be\label{estionep2}
\frac{ ( | \om -1 | n^2 )^{(p-1)\slash (d-1)}}{\gamma_\om}  \leq C_{11}  
\ee
and $ \om > 1 $ if $b >0 $ (resp. $ \om <1 $ if $ b  < 0 $),
equation (\ref{eq:main}) possesses at least one pair of even
periodic (in time) classical $C^2 $ 
solutions of minimal period $2\pi / (n\om)$.
\end{theorem} 

\begin{pf}
In this case 
the reduced functional $ \Phi_{\e} $ can be written in the form 
(\ref{absfunc}) with $\mu=\e $, $ G(v) := $ $(b/ d+1) $ 
$ \int_{\Omega} v^{d+1} $,
where $b := f^{(d)} ( 0 ) \slash d! $, and 
by lemma \ref{appro} ( since $ \int_\Omega v^{d+2} = 0) $
$$
| D R (v)[v]| \leq C \Big( |v|_\om^{d+3} + \frac{\e}{\gamma}
|v|_\om^{2p} + \frac{|v|_\om^{3p-1}}{\gamma} \Big).
$$
As before $G \geq 0 $ and $ G({\cal L}_n v) = G( v ) $.
Applying the same arguments as in the previous subsection, with
$ q = d $, $ r_0 = (1 \slash C_0 ) ( |\e | n^2 \slash m )^{1/d-1} $, 
and 
$$
\alpha = O \Big( 
{r_0}^2+ \frac{| \e |}{\gamma} {r_0}^{2p-1-d} + \frac{{r_0}^{3p-2-d}}{\gamma}
\Big),
$$
we obtain that condition (\ref{estionep2}) ensures that 
$B_{r_0} \subset {\cal D} $ and that conditions (\ref{estionep})
hold, which implies the existence of
a pair of critical points of $\Phi_{\e,n}$ (observe that 
$2p-1-d \geq 2$).
\end{pf}

When $(N1)$
is not satisfied the situation is  more delicate.
We shall first consider the case in which $(N2)$ holds.
We can write, since $ \Pi_W v^p = v^p $,  
$$
\frac{1}{2} \int_\Om f(v) L^{-1} \Pi_W f(v) = - G(v)+ r(v),
$$ 
with 
\be\label{Geven}
G(v) = - \frac{a^2}{2} \int_\Om v^p L^{-1} v^p \qquad {\rm and} \qquad
Dr(v)[v]= O(|v|_\om^{2p+2})
\ee
(recall that $f^{(p+1)}(0) = 0 $ by $(N2)$).
Hence, by lemma \ref{appro}, we have
$ - \Phi_\e (v) = (-\e \slash 2 )||v||^2 - G(v) + R(v),$
with
\be\label{DReven}
DR(v)[v]= O \Big( |v|_\om^{2p+2} + \frac{\e}{\gamma} |v|_\om^{2p} +
\frac{|v|_\om^{3p-1}}{\gamma} \Big).
\ee
$ D G, DR : V \to V^* $ are still compact operators  
and now $ G $ is homogeneous of degree $ q + 1 := 2 p $.
\\[1mm]
\indent
A difficulty arises for proving that $ G $ satisfies
 (\ref{smallGn}). 
The following lemmas will enable to 
get a suitable expression of $ G ( v )$.

\begin{lemma} \label{lemexpl} 
Let $ w(t,x) \in W $  be of the form 
$ w (t,x) = m(t+x,t-x) $ for some 
$ m : {\bf R}^2 \to {\bf R}$ $2\pi$-periodic
with respect to both variables. 
Then there exists $\wtilde{m} $, $ a $ $2-\pi$ periodic  
such that $ m $ can be decomposed as 
\be\label{decm}
m ( s_1, s_2 ) = \wtilde{m} (s_1, s_2) + a(s_1)+a(s_2)+ \alpha,
\ee
where $\alpha \in {\bf R}$ 
($\alpha= \la m \ra $), $ \la a \ra =0$, $ \la \wtilde{m} \ra_{s_1} (s_2)=0$,
$\la \wtilde{m} \ra_{s_2} (s_1)=0$. Moreover
\begin{eqnarray*} \label{calculpenible}
\int_0^{2\pi} \int_0^{\pi} L^{-1} (w) w \ dt \ dx
&=& - \frac{1}{2} \int_0^{2\pi} \int_0^{2\pi} M(s_1,s_2)
\wtilde{m}(s_1,s_2) \ ds_1 \ ds_2 +
2\pi \int_0^{2\pi}  M(s,s)
a(s) \ ds \\
&+& 2\pi \alpha \int_0^{2\pi}  M(s,s) \ ds
- 8\pi \int_0^{2\pi}  A(s)^2
 \ ds - \frac{\alpha^2 \pi^4}{6},
\end{eqnarray*}
where $M$ and $A$ are the $2\pi$-periodic 
functions defined by
\be\label{eq:AM}
A'(s)=\frac{1}{4}a(s), \ \  \la A \ra =0, \ \
\partial_{s_1} \partial_{s_2}  M (s_1,s_2) = 
\frac{1}{4}\wtilde{m}(s_1,s_2), \ \
\la M \ra_{s_1} (s_2) =0, \ \
\la M \ra_{s_2} (s_1) =0.
\ee
\end{lemma}

\begin{pf} 
The proof is in the Appendix.
\end{pf}

\begin{lemma}\label{cor:nl}
Let $ m : {\bf R}^2 \to {\bf R}$ be $2\pi$-periodic
with respect to both variables and such that
$w(t,x)=m(t+x,t-x) \in W$. For $ n\in {\bf N}$  define 
$({\cal L}_n w) (t,x) := m(n(t+x),n(t-x))$. Then 
\be\label{rescaln}
\int_0^{2\pi} \int_0^{\pi} L^{-1} ({\cal L}_n w) {\cal L}_n w \ dx \
dt=-\frac{\pi^4}{6} \alpha^2 + \frac{1}{n^2} 
\Big( \int_0^{2\pi} \int_0^{\pi} w L^{-1} w \ dx \ dt +
\frac{\pi^4}{6} \alpha^2 \Big).
\ee
\end{lemma}
\begin{pf}
We have $ {\cal L}_n w (t,x) = m_n (t+x, t-x) $, with
$ m_n(s_1,s_2) := m (ns_1,ns_2)$. Using the decomposition given in the
latter lemma, we can write $ m_n(s_1,s_2) = \wtilde{m}_n ( s_1, s_2 )+
a_n (s_1)+a_n(s_2) + \alpha$, where $ \wtilde{m}_n :=m( n \ \cdot )$ and 
$a_n := a(n \ \cdot )$.
Therefore (using the abbreviation $ w_n := {\cal L}_n w$)
\begin{eqnarray*}
\int_0^{2\pi} \int_0^{\pi} L^{-1} (w_n) w_n \ dx \ dt
&=& -\frac{1}{2} \int_0^{2\pi} \int_0^{2\pi} M_n(s_1,s_2)
\wtilde{m}_n(s_1,s_2) \ ds_1 \ ds_2 +
2\pi \int_0^{2\pi}  M_n(s,s)
a_n(s) \ ds \\
&+& 2\pi \alpha \int_0^{2\pi}  M_n(s,s) \ ds
- 8\pi \int_0^{2\pi}  A_n(s)^2
 \ ds - \frac{\alpha^2 \pi^4}{6},
\end{eqnarray*} 
where $M_n(s_1,s_2) := M(ns_1,ns_2)/n^2$ and  $A_n(s) :=A(ns)/n$. 
 (\ref{rescaln}) follows straightforwardly.
\end{pf}

Lemma \ref{lemexpl} allows to obtain a handable expression for 
$ G(v) = - (a^2 \slash 2)  \int v^p L^{-1}v^p $ when $p=2$. 
For $ v (t, x) = \eta (t+x) - \eta (t-x) 
\in V $ we have $ v^2 (t,x) = m (t+x, t-x) \in W $
with $ m(s_1, s_2) = ( \eta (s_1) - \eta (s_2) )^2 $. 
We can then decompose $ m $ as in (\ref{decm}) with
$ \alpha = 2 \la \eta^2 \ra $,  
$ a ( s ) = \eta^2 (s) -  \la \eta^2 \ra $ and 
$ \wtilde{m}(s_1, s_2 ) = - 2 \eta (s_1) \eta ( s_2 ) $.
Let us call 
$ P_1 $, resp. $ P_2 $, the primitive of $ \eta $, 
resp. $ \eta^2 - \la \eta^2 \ra $, of zero mean-value. 
The functions $ A, M $ defined in (\ref{eq:AM})
are $ A(s) := P_2 (s) \slash 4 $ and 
$ M(s_1, s_2) = - P_1 ( s_1 ) P_1 ( s_2 ) \slash 2 $.    
Hence, by lemma \ref{lemexpl} we obtain the following expression
for $ G(v) = - (a^2 \slash 2) \int v^2 L^{-1}v^2 $ :
\be\label{p=2}
G(v) = \frac{a^2}{2} \Big[
\pi \int_0^{2\pi} P_1(s)^2 \Big( \eta(s)^2 + \la \eta^2 \ra \Big)
\ ds + \frac{\pi}{2} \int_0^{2\pi} P_2(s)^2 \ ds +
\frac{2\pi^4}{3} \la \eta^2 \ra^2 \Big]. 
\ee
From (\ref{p=2}) we see that 
\be\label{pop=2} 
\forall v \in V, \ G (v ) > 0 \qquad {\rm and } \qquad
0 \leq G ( {\cal L}_n v ) \leq G ( v ), \ \forall n \geq 2.
\ee 
In general we can prove that
\begin{lemma}\label{Gpos}
For any even $ p $,
$ G(v) := - (a^2 \slash 2 )  \int_\Omega v^p L^{-1}v^p $ is non negative.
\end{lemma}

\begin{pf}
Writing $ v^p(t,x)=m(t+x,t-x) $ with 
$ m(s_1,s_2) = (\eta (s_1)-\eta(s_2))^p $, we claim 
that $L^{-1}(v^{p}) \in $ $ -( 1/ 8) M(t+x,t-x)+ V $, where
\be \label{positivity}
M(s_1,s_2) := \int_{R_{s_1,s_2}} m (\xi, \nu) \ d\xi \ d\nu,
\ee
and
$ R_{s_1,s_2} := $ $\{ (\xi,\nu) \in {\bf R}^2 \ | \ s_1 \leq \xi \leq$ 
$s_2+2\pi, \ s_2 \leq \nu \leq s_1 \}.$ In fact 
$$
\partial_{s_1} M(s_1, s_2 ) = \int_{s_1}^{s_2+2\pi} m(\xi,s_1) \ d\xi
-\int_{s_2}^{s_1} m(s_1,\nu) \ d\nu.
$$
Hence, since $ m(s_1,s_2) = (\eta (s_1)-\eta(s_2))^p$ with even $p$,
and $\eta$ is $2\pi$-periodic,
$$
\partial_{s_2} \partial_{s_1} M = m(s_2+2\pi,s_1) +m(s_1,s_2)=
2 m(s_1,s_2).
$$
Moreover, by the definition of $M$, 
$M(s_1,s_1)=M(s_1,s_1-2\pi)=0$ and 
$ M (s_1 + 2 \pi, s_2 + 2 \pi)= M(s_1, s_2) $ for all 
$(s_1, s_2) \in {\bf R}^2 $,  
so that $(t,x) \mapsto M(t+x,t-x) \in X$.
Hence $L [M(t+x,t-x)]= -4 (\partial_{s_2} \partial_{s_1} M) 
(t+x,t-x)=- 8 v^p (t,x)$ and our claim follows.
From (\ref{positivity}) we see that 
$M$ is non negative and as a result 
$ - a^2 \int L^{-1} (v^p) v^p = (a^2 \slash 8 ) \int M v^p \geq 0$.
\end{pf}

We do not know whether (\ref{pop=2}) still  holds true for 
$ G(v) = - (a^2 \slash 2) \int v^p L^{-1} v^p $,
$ p \geq 4 $. In any case, we can deduce from lemma \ref{Gpos},
 that $ G_n (v) := G ({\cal L}_n v ) $ satisfies (\ref{smallGn})
for $ n \geq 2 $. Indeed, by (\ref{rescaln}),  
\be\label{Gn}
G_n(v) = G( {\cal L}_n v ) = \frac{G(v)}{n^2} + 
\frac{ a^2 \pi^4 \alpha^2 }{12}
\Big( 1 - \frac{1}{n^2} \Big ),
\ee
with $ \alpha := (1\slash 2 \pi^2 ) \int_\Omega v^p $.
Hence, since $G(v) \geq 0$, $ \forall n, m \geq 2 $, 
for all $v\in V \backslash \{ 0 \}$,
\be\label{rescalepn}
G_n ({\cal L}_m v) =
G( {\cal L}_{mn} v) = \frac{G(v)}{(nm)^2} + \frac{a^2 \pi^4 \alpha^2 }{12}
\Big( 1 - \frac{1}{(nm)^2} \Big )\leq \frac{G(v)}{n^2}
+\frac{a^2 \pi^4 \alpha^2}{12} <  m^{2p} G_n( v). 
\ee 
We may write 
for $ v(t,x) = \eta ( t + x ) - \eta ( t - x ) \in V $
$$
-\Phi_{\e,n}(v) = \frac{- \e n^2}{2}||v||^2- G_n (v )+R_n(v),
$$ 
where $ G_n (v) = G ( {\cal L}_n v )$, $ R_n(v)= R({\cal L}_n v ) $ 
satisfies, by (\ref{DReven}),     
\be \label{restepair}
DR_n(v)[v]=DR({\cal L}_n v)[{\cal L}_n v]
= O \Big( ||v||^{2p+2} + \frac{\e}{\gamma} ||v||^{2p} +
\frac{||v||^{3p-1}}{\gamma}  
\Big).
\ee
We shall consider frequencies $ \om <1 $ so that $ -\e >0 $.
By (\ref{rescalepn}), 
for $ n \geq 2 $, $ G_n (v) $ satisfies (\ref{smallGn}). 
For $(|\e | n^2)^{1/2} \slash \gamma $ small enough we shall apply      
Proposition \ref{mainlem} and lemma \ref{miniper} 
to $ -\Phi_{\e,n} : B_{r_0} \to {\bf R} $ 
as before, with $ q = 2p-1 $, $\mu = -\e n^2 $, 
$$ 
{r_0}= \frac{1}{C_0} \Big( \frac{|\e | n^2}{m( G_n )} \Big)^{1/(2p-2)} \qquad
{\rm  and} \qquad  \alpha = O \Big( {r_0}^2+ \frac{| \e |}{\gamma}
+\frac{{r_0}^{p-1}}{\gamma} 
\Big).
$$
Note that $m(G_n) \to m(\ov{G})$ as $n\to +\infty$, where
$\ov{G} (v):= (a^2 \pi^4/ 12) \la v^p \ra^2 = (a^2 /48)
(\int v^p )^2$, and the sequence $(m(G_n))$ is bounded 
from below by a positive constant.
We derive the following Theorem, where the condition 
$r_0^{p-1}\slash \gamma $ small, i.e. 
$(| \e | n^2 )^{1 \slash 2} \slash \gamma$ small, is required.

\begin{theorem} \label{theven}
Assume that  $f$ satisfies $(F1)$ and $(N2)$. 
Then there is a constant $C_{12}$ depending only
on $f$ such that, $ \forall \om \in {\cal W}$, $ \om < 1 $,  
$ \forall n \geq 2 $ such that 
\be\label{estionep3}
\frac{(|\om -1| n^2)^{1 \slash 2}}{\gamma_\om} \leq C_{12}
\ee
equation (\ref{eq:main}) possesses at least one pair of even
periodic (in time) classical $ C^2 $ 
solutions of minimal period $2\pi / (n \om )$.
If $ p = 2 $ the existence result  holds true 
for  $ n =1 $ as well.
\end{theorem} 

\begin{remark}
We think that the restriction $ n \geq 2 $  
in the case $ p \neq 2 $ even is of purely technical nature.
\end{remark}

Finally we consider the case in which $(N3)$ holds.
We first assume that $b<0$. We have
$$
-\Phi_{\e}(v) = \frac{-\e }{2}||v||^2- G(v)+R(v),
$$
where
\be\label{GN3}
G(v):= - \frac{b}{2p} \int_{\Omega} v^{2p} - 
\frac{a^2}{2} \int_{\Omega} v^p L^{-1}v^p >0
\ee
for all $v\neq 0$ and $ R(v) $ satisfies  estimate 
(\ref{DReven}). Moreover, defining $ G_n(v) := G( {\cal L}_n v) $,
we  still have $G_n ({\cal L}_m v)
< m^{2p} G_n (v)$ for all $ n, m \geq 2 $ and $ v\neq 0 $. For
$\e <0$ and $|\e |n^2 \slash \gamma $ small enough we 
obtain the same existence result as in Theorem \ref{theven}.
\\[3mm] 
\indent Now assume that $ b > 0 $. 
As a consequence of lemma \ref{cor:nl} we have
\be\label{smallresca}
\int_{\Omega} {\cal L}_n w L^{-1} ({\cal L}_n w) =-\frac{\pi^4}{6} \la m
\ra^2 + O \Big( \frac{|m|^2_{L^2} + |w|^2_{L^2} }{n^2} \Big) 
= - \frac{\pi^4}{6} \la w
\ra^2 + O \Big( \frac{|w|^2_{L^2}}{n^2} \Big).
\ee
By (\ref{smallresca}), we may write 
for $ v(t,x) = \eta ( t + x ) - \eta ( t - x ) \in V $
$$
\Phi_{\e,n}(v) = \frac{\e n^2}{2}||v||^2-\wtilde{G}(v)+R_n(v),
$$ 
with 
$$
\wtilde{G}(v)=\frac{b}{2p} \int_\Om v^{2p} - \frac{a^2}{48} 
\Big(\int_\Om v^p  \Big)^2 = \frac{b}{4p} \int_{{\bf T}^2}
\Big( \eta(s_1)-\eta(s_2) \Big)^{2p} - \frac{a^2}{192} \Big( \int_{{\bf T}^2} 
(\eta(s_1)-\eta(s_2))^p \Big)^2 
$$
and
\be \label{restepairn}
DR_n(v)[v]=DR({\cal L}_n v)[{\cal L}_n v]+ 
O \Big( \frac{|v|_{\om,n}^{2p}}{n^2} 
\Big) = O \Big( ||v||^{2p+2} + \frac{\e}{\gamma} ||v||^{2p} +
\frac{||v||^{3p-1}}{\gamma}+\frac{||v||^{2p}}{n^2} \Big).
\ee 
We are no longer able to prove that the condition of lemma 
\ref{smallGn} is satisfied by $ G $ defined in (\ref{GN3})
or by $ - G $. That is why we introduce $\wtilde{G}$.
Clearly $ { \wtilde G} ({\cal L}_n v ) = {\wtilde G} ( v ) $ for all
$v \in V $ and $  n \in {\bf N}\backslash \{0\} $.  
We  have to specify the sign of $\wtilde{G}(v)$. 
We need to introduce
$$
\kappa(p):= \sup_{v \in V\backslash \{ 0\}} \Big( \int_\Om v^p
\Big)^2 \Big( \int_\Om  v^{2p} \Big)^{-1} =
\sup_{\eta \in H^1({\bf T}), \eta \ odd} \Big( \int_{{\bf T}^2} 
(\eta(s_1)-\eta(s_2))^p
\Big)^2 \Big(2 \int_{{\bf T}^2}  (\eta(s_1)-\eta(s_2))^{2p} \Big)^{-1}.
$$
We can prove that $ \kappa(p) = \pi^2 $. 
Indeed, using that $ \eta $ is odd, we get
$$
 \int_{\T \times \T}
(\eta(s_1)-\eta(s_2))^p \ ds_1 \ ds_2 = \sum_{k=0}^{p/2}
C_p^{2k}  \int_{\T} \eta^{2k} \ ds  \int_{\T} \eta^{p-2k} 
\ ds \leq \sum_{k=0}^{p/2}
C_p^{2k} \times 2\pi \int_{\T} \eta^p \ ds,
$$ 
by H\"older inequality. Hence 
$$
 \int_{\T \times \T}
(\eta(s_1)-\eta(s_2))^p \ ds_1 \ ds_2 \leq 2^p \pi \int_{\T} \eta^p \ ds.
$$
Now, by Cauchy-Schwarz inequality, for all $ 0 \leq k \leq p$,
$$
\Big( \int_{\T} \eta^p \Big)^2 = \Big( \int_{\T} |\eta|^k
|\eta|^{p-k} \Big)^2 \leq  \int_{\T} \eta^{2k} 
 \int_{\T} \eta^{2p-2k}.
 $$
Hence 
$$
 \int_{\T \times \T}
(\eta(s_1)-\eta(s_2))^{2p} \ ds_1 \ ds_2 =
\sum_{k=0}^{p}
C_{2p}^{2k}  \int_{\T} \eta^{2k} \ ds  \int_{\T} \eta^{2p-2k} 
\ ds \geq 2^{2p-1} \Big( \int_{\T} \eta^p \Big)^2.
$$
As a result we obtain $\kappa(p) \leq \pi^2$. Choosing for 
$\eta$ $H^1$-approximations of the $2\pi$-periodic function 
that takes value $-1$ on $(-\pi,0]$ and value $1$ on $(0,\pi]$,
we obtain the converse inequality. Hence $\kappa (p)=\pi^2$.
Note also that
$$
\inf_{v \in V\backslash \{ 0\}} \Big( \int_\Om v^p
\Big)^2 \Big( \int_\Om  v^{2p} \Big)^{-1} =0.
$$
In fact, for all $\delta>0$ there is an odd smooth
$2\pi$-periodic  $\eta_\delta$ such that $\int_0^{2\pi} 
\eta_\d^{2p} =1$ and $\int_0^{2\pi} \eta_\d^{2q} \leq \delta$ for
$1\leq q < p$. One can see easily that for 
$v(t,x)=\eta_\delta (t+x)-\eta_\delta (t-x)$, 
$( \int_\Om v^p
)^2 ( \int_\Om  v^{2p} )^{-1} \leq \beta(\delta)$ with 
$\lim_{\delta\to 0} \beta(\delta)=0$.

If $ b\geq p \kappa(p) a^2/24=p \pi^2 a^2/ 24$ then for all $v\in V$,
$\wtilde{G}(v) \geq 0$ and $G \not\equiv 0$. In this case we
obtain an existence result for $\om >1$.

If $0<b<p\kappa(p)a^2/24 = p \pi^2 a^2/ 24 $ then there is $v_1 \in V$ such
that $\wtilde{G}(v_1)<0$ and there is $v_2 \in V$ such that
$\wtilde{G}(v_2)>0$. In this case we obtain an existence result
both for $\om <1$ and $\om>1$. 

In both cases  we  apply 
the same arguments as before, for $(|\e |n^2 )^{1/2}\slash \gamma $ small
and $n$ large enough,   
with $ q = 2p - 1 $, $\mu = |\e | n^2 $, 
$$ 
{r_0}= \frac{1}{C_0} \Big( \frac{| \e | n^2}{m(\pm \wtilde{G})} 
\Big)^{1/(2p-2)} \qquad
{\rm  and} \qquad  \alpha = O \Big( {r_0}^2+ \frac{\e}{\gamma}
+\frac{{r_0}^{p-1}}{\gamma} + \frac{1}{n^2} \Big).
$$
$m(-\wtilde{G})$ appears only in the second case when 
$\om < 1$.
As a conclusion we have

\begin{theorem}\label{thm:main3}
Assume that $f$ satisfies $(F1)$ and $(N3)$.
Define   $J \subset {\bf R} $ as $ J = ( 0, 1 ) $ if $ b < 0 $,
$J=(1,2)$ if $ b \geq p \pi^2 a^2/ 24$ and $J=(0,1)\cup (1,2)$ if
$0<b<p \pi^2 a^2/24$.
Then there is a constant $ C_{13} $ and an
integer $ N_0 \geq 1$ depending only
on $ f $ (with $N_0=2$ if $b<0$ and $p\geq 4$, 
$N_0=1$ if $ b < 0$ and $ p = 2$)
such that, $ \forall \om \in {\cal W}\cap J $,   
$\forall n \geq N_0 $ such that 
\be\label{estionep4}
\frac{(| \om - 1 | n^2)^{1 \slash 2}}{\gamma_\om} \leq C_{13}
\ee
equation (\ref{eq:main}) possesses at least one pair of even
periodic (in time) classical $ C^2 $ 
solutions of minimal period $2\pi / (n \om )$.
\end{theorem}

In the next remarks we specify the norms and the 
energies of the solutions $ u_n $ obtained
in Theorems \ref{thm:main1}, \ref{thm:main2}, \ref{theven},  
\ref{thm:main3}.

\begin{remark}\label{rem:red}
It is readily seen by the choice of $r_0 $ in the
different cases and the fact that 
$ | {\cal L}_n v|_\om =$ $ O( ||v|| ) = $ $O(r_0) $, for 
$ |\e | n^2 < 1 $, that the critical points $v_{\e,n} $ 
of $\Phi_{\e}  $ that we obtain satisfy 
$ |v_{\e,n}|_\om = O((|\e| n^2)^{1\slash (q-1)})$ with
$q = p $ if $ p $ is odd and $q= d$ in case $(N1)$,  
$q= 2p-1$ in cases $(N2)$-$(N3)$, if $p$ is even. 
\end{remark}

\begin{remark}\label{rem:nen}
Let $ u_{\om,n}:= u_n $ be the sequence of solutions ($u_{\om,n}$ has minimal 
period $2\pi/n \om $) obtained from the above Theorems. Let 
${\cal E}_{\om,n} := 
\int_0^\pi (u_{\om,n})_t^2/2 + (u_{\om,n})_x^2/2 + F(u) \
dx$ be the energy of $u_{\om,n}$. It is possible to prove, using Proposition
(\ref{mainlem}), the
following estimates:
\be\label{estino}
||u_{\om,n}||_{H^1} = n \Big(\frac{|\om-1| n^2}{ m_n (q+1)} 
\Big )^{1/(q-1)} \Big( 1+ g_1 \Big( 
\frac{(|\om-1| n^2)^{(p-1)\slash (q-1)}}{\gamma_\om} \Big) \Big)
\ee
$$  
{\cal E}_{\om,n}=\frac{n^2}{2\pi} \Big(\frac{|\om-1| n^2}{m_n (q+1)} 
\Big )^{2/(q-1)} \Big( 1+ g_2 \Big( 
\frac{(|\om-1| n^2)^{(p-1)\slash (q-1)}}{\gamma_\om} \Big) \Big)
$$ 
where $q=p$ if $p$ is odd, $q=d$ if $p$ is even and $d<2p-1$ 
is the smallest order of the odd terms, $q=2p-1$ if $p$ is even
and there is no odd term of order $<2p-1$ and $\lim_{s\to 0}
g_i(s)=0$; $m_n := m(G_n) =m$ if $p$ is odd and $m_n=m+O(1/n^2)$ if 
$p$ is even, where $m>0$ depends only on $f$. 
In the  case when the smallest order of the odd terms
is $2p-1$ and $0<b<p\pi^2 a^2/24$, 
we have two sequences of solutions which correspond to
two (possibly) different values of $m$.  
\end{remark}

We finally discuss the optimality of the estimate 
on the number of solutions that we find. 

\begin{remark}\label{rem:opt}
By conditions
(\ref{estionep1}), (\ref{estionep2}), 
(\ref{estionep3}), (\ref{estionep4})
in Theorems \ref{thm:main1}, \ref{thm:main2}, \ref{theven},  
\ref{thm:main3} we can find 
$ N_\om \approx O(\gamma_\om^\tau \slash \sqrt{|\om - 1|})$ with 
$1 \leq \tau \leq 2 $, periodic solutions
of (\ref{eq:main}) with period $ 2 \pi \slash \om $.  
They have decreasing minimal periods $ 2 \pi \slash \om, \pi \slash \om, 
\ldots, 2 \pi \slash N_\om \om $   
and increasing energies in $ n $, see remark \ref{rem:nen}.
We expect the estimate on the number $ N_\om $ of solutions found, with 
small $| \cdot |_\om^{p-1} \slash \gamma_\om $, to be optimal since 
we could have performed a finite
dimensional Lyapunov-Schmidt reduction $ u := v_N + w $
with $ v_N := \sum_{j=1}^N \xi_j \cos jt \sin jx $ and 
$ N^2 |\om - 1 | = O( \gamma_\om )$. Indeed, 
in lemma \ref{lem:red}, the ``small denominators'' are
$ | \om^2 l^2 - j^2 | \geq \gamma_\om $, $ \forall l \neq j $, 
and for $ j = l \geq N $ are $ | \om^2 l^2 - j^2 | =
| \om^2 - 1 | j^2 \geq ( |\om - 1| \slash 2 ) N^2 = O( \gamma_\om ) $.
In this way we reduce to study a 
even functional like $ \Phi_\e $ on a finite dimensional 
space of dimension $ N \approx O(\gamma_\om \slash \sqrt{|\om - 1|})$. 
In general the best estimate regarding the 
number of its critical points is $ N $.
\end{remark}

\section{Regularity}\label{sec:reg}

In order to complete the proof of Theorems \ref{thm:main1},
\ref{thm:main2}, \ref{theven}, \ref{thm:main3} 
we now prove that the weak solutions $ u_n $  
actually are classical $ C^2 ( \Omega )$ solutions of 
equation (\ref{eq:main}). 
\\[1mm]
\indent
Define the norm
$$ 
| u |_{\om,3} := |u|_\infty + \sqrt{|\e |} \  || u ||_{H^1} +
|\e | \ || u ||_{H^2} + |\e |^{3/2} \ || u ||_{H^3}.
$$
It is easy to see, arguing as in lemma \ref{composition}, and 
using the Sobolev continuous embedding 
$H^3 (\Om) \subset W^{1,\infty}(\Omega) $ 
that 
\be\label{nemint3}
|f(u)|_{\om,3} \leq C | u |_{\om,3}^p . 
\ee
The following lemma holds:

\begin{lemma}\label{regul}
If $ u \in X \cap H^3(\Omega) $  then 
$
w:= L_\om^{-1} \Pi_W u := \sum_{l \geq 0, j \geq 1, j \neq l} 
w_{lj} \cos (lt) \sin (jx) $, where 
$ w_{lj} := $ $ u_{lj} \slash (\om^2 l^2 - j^2) $, 
satisfies $ \sum_{l \geq 0, j \geq 1, j \neq l} 
|w_{lj}|^2 (l^6 + j^6 ) \leq (C_{14} \slash \gamma )^2 ||u||_{H^3}^2 $.
In particular $w\in H^3(\Omega)$, $||w||_{H^3}=
O(||u||_{H^3})/\gamma$ and  $ w_{xx}(t,0)= w_{xx}(t,\pi) = 0 $.
Moreover $ w \in   C^2 (\ov{\Omega} ) $.
 
\end{lemma}

\begin{pf}
We point out that 
$ u \in H^3 (\Omega ) $ does not imply that 
$ \sum_{l \geq 0, j \geq 1}|u_{lj}|^2  j^6  $
is finite (it is right only if also $u_{xx}(t,0) = u_{xx}(t,\pi) = 0$;
in general we have only  $||u||_{H^3}^2 \leq C \sum_{l \geq 0, j \geq 1} 
|u_{lj}|^2 (l^6 + j^6 )$).
However it is true that  
$ \sum_{ l \geq 0, j \geq 1, j \neq l } 
|u_{lj}|^2 l^6 = $ $ O( | \partial_{ttt} u |_{L^2}^2)
=O(||u||^2_{H^3}) $. 
Now
\be\label{eq:2}
\sum_{ l \geq 0, j \geq 1, j \neq l} w_{lj}^2 (l^6+ j^6) = 
\sum_{ l \geq 0, j \geq 1, j \neq l} 
\frac{|u_{lj}|^2}{(\om^2 l^2 - j^2)^2}(l^6+ j^6) = S_1 + S_2
\ee
where 
$$
S_1 := \sum_{j \geq 2 \om l } 
\frac{|u_{lj}|^2}{(\om^2 l^2 - j^2)^2}(l^6+ j^6)
 \qquad {\rm and} \qquad  
S_2 := \sum_{j < 2\om l, j\neq l } 
\frac{|u_{lj}|^2}{(\om^2 l^2 - j^2)^2}(l^6+ j^6). 
$$
If $ j \geq 2 \om l $ then $ |\om^2 l^2 - j^2| \geq 3 j^2 \slash 4 $,
hence
\be\label{eq:3}
S_1 \leq \sum_{j \geq 2 \om l } 
\frac{2 u_{lj}^2}{ j^4} \Big( j^6 + \frac{j^6}{(2\om)^2} \Big) \leq 
\sum_{j \geq 2 \om l } 4|u_{lj}|^2 j^2 =O( | u_{x} |^2_{L^2})
=O( ||u||^2_{H^1}).
\ee
Since $| \om^2 l^2 - j^2 | \geq \gamma /2 $ for all $j \neq l $ and 
$ \om \leq 3 \slash 2 $ we have 
\be\label{eq:4}
S_2 \leq \sum_{1\leq j < 2 \om l } 
\frac{4|u_{lj}|^2}{\gamma^2}(l^6+ j^6) \leq \sum_{1\leq j < 2 \om l } 
\frac{4|u_{lj}|^2}{\gamma^2 }(l^6+ (3l)^6) \leq
\frac{C}{\gamma^2} |u_{ttt}|_{L^2}^2
\leq \frac{C'}{\gamma^2} ||u||_{H^3}^2.
\ee
By  (\ref{eq:2}),
(\ref{eq:3}) and (\ref{eq:4}) we get 
 $ \sum_{l \geq 0, j \geq 1, j \neq l} 
|w_{lj}|^2 (l^6 + j^6 ) \leq (C \slash \gamma )^2 ||u||_{H^3}^2 $.

Now we prove that 
\be\label{der2}
\sum |w_{lj}|(l^2+j^2) < \infty,
\ee
which implies that $w\in C^2(\ov{\Omega})$. We have
\be\label{s1s2}
\sum_{ l \geq 0, j \geq 1, j \neq l} |w_{lj}| (l^2+j^2) = 
\sum_{ l \geq 0, j \geq 1, j \neq l} 
\frac{|u_{lj}|}{|\om^2 l^2 - j^2|}  (l^2+j^2) = S'_1 + S'_2
\ee
where 
$$
S'_1 := \sum_{j \geq 2 \om l } 
\frac{|u_{lj}|}{|\om^2 l^2 - j^2|}(l^2+ j^2) \qquad {\rm and} \qquad  
S_2' := \sum_{j < 2\om l, j\neq l } 
\frac{|u_{lj}|}{|\om^2 l^2 - j^2|} (l^2+j^2).
$$
For $ j \geq 2 \om l $, $ |\om^2 l^2 - j^2| \geq 3 j^2 \slash 4 $
and $ l \leq j \slash 2 \om \leq j $,  hence
\be\label{S1'}
S'_1 
\leq C \sum_{j \geq 2 \om l } |u_{lj}| \leq C
\Big( \sum_{l\geq 0, j \geq 1} |u_{lj}|^2 (l^2 + j^2)^2 \Big)^{1 \slash 2}
\Big( \sum_{l\geq 0, j \geq 1} \frac{1}{(l^2 + j^2)^2} \Big)^{1 \slash 2} 
\leq C' ||u||_{H^2}.
\ee
Note that $\sum_{l\geq 0, j \geq 1} |u_{lj}|^2 (l^2 + j^2)^2 
= O( ||u||_{H^2}^2 )$. We claim that
\be\label{S2'}
S_2' \leq
 \sum_{j < 2\om l, j\neq l } C 
\frac{|u_{lj}|l^2}{|\om^2 l^2 - j^2|} < + \infty.
\ee
Indeed we know that
\be \label{auxil}
\sum_{j < 2\om l } 
(u_{lj} l^2)^2 l^2  < + \infty.
\ee
As in the proof of (\ref{eq:S}) in lemma \ref{lem:operL}, 
(\ref{auxil}) implies that $ S'_2 $ is finite. 
By and (\ref{s1s2}), (\ref{S1'}) and (\ref{S2'}) we get (\ref{der2}).
The lemma is proved.
\end{pf}

By lemma \ref{regul}
$ || L_\om^{-1} \Pi_W ||_{\om,3} \leq C \slash \gamma $. Hence, 
by (\ref{nemint3}), 
the map $ { \cal G }_v : W \to W $, defined 
in the proof of lemma \ref{lem:red} by
$ { \cal G}_v (w) :=   L_\om^{-1} \Pi_W f(v + w) $,
satisfies (if $v,w \in H^3(\Omega)$)
\be\label{smallg}
|{ \cal G}_v (w)|_{\om,3} \leq \frac{\wtilde{C}}{\gamma} 
\Big( | v |_{\om,3}^p + | w |_{\om,3}^p \Big),
\ee
where $\wtilde{C}$ is some positive constant.
Define $\wtilde{\d}:= [\gamma/ (2\wtilde{C})]^{1/(p-1)}$ and assume that
$|v|_{\om,3} \leq \wtilde{\d}$. Then by (\ref{smallg}),
${\cal G}_v (\wtilde{B}_{\wtilde{\d}}) \subset 
\wtilde{B}_{\wtilde{\d}}$, where
$$
\wtilde{B}_{\wtilde{\d}}:= \Big\{ w \in W \cap H^3(\Omega) \
\Big|  \  |w|_{\om,3} \leq \wtilde{\d}  \Big\}.
$$
As a result, if $|v|_{\om,3} \leq \wtilde{\d}$ then 
$w_k := {\cal G}_v^k (0) \in \wtilde{B}_{\wtilde{\d}}$
for all $k\geq 1$. It follows that $(w_k)$ converges weakly
in $H^3(\Omega)$ to some $\wtilde{w} \in \wtilde{B}_{\wtilde{\d}}$.

Now, recall that $ w(v) $ is the unique fixed point of 
$ { \cal G }_v $ in a small ball around $0$ for the norm 
$| \cdot |_\om $ in $W$ and that, by the Contraction Mapping Principle,
 $w(v) $ is the $ H^1 $ limit of $ (w_k)  $.
Hence $w(v)= \wtilde{w}$ and $w(v) \in H^3(\Omega)$. 
By formula (\ref{nemint3}) and lemma \ref{regul} we derive from
the above considerations that if $|v|_{\om,3} \leq \wtilde{\d}$
then $ w(v) \in C^2 (\ov{\Omega}) $, and $v+w(v) \in  C^2 (\ov{\Omega})$
(because $V\cap H^3(\Omega) \subset  C^2 (\ov{\Omega}))$.

There remains to check that when condition
(\ref{estionep1}) or (\ref{estionep2}) or
(\ref{estionep3}) or (\ref{estionep4}) is satisfied
(according to the different cases)
the critical points $v_{\e ,n}$ of $\Phi_{\e}$ that we have obtained
in section \ref{secappl} satisfy $|v_{\e,n}|_{\om,3}^{p-1} \leq 
\gamma/(2\wtilde{C})$.
Since the critical points $ v := v_{\e,n} $ of 
$ \Phi_\e $, obtained in Theorems \ref{thm:main1},
\ref{thm:main2}, \ref{theven}, \ref{thm:main3} 
satisfy $ v_{tt} = v_{xx} $,  
$ 2 \e v_{tt} = \Pi_V f(v + w(v)) $,
$$
||v||_{H^3} = O\Big( \frac{1}{\e} || \Pi_V f(v + w(v))|| \Big) =
O\Big( \frac{1}{\e}|| f(v + w(v))||  \Big), \quad
||v||_{H^2}  =
O\Big( \frac{1}{\e}| f(v + w(v))|_{L^2} \Big). 
$$
Therefore by lemmas \ref{composition} and \ref{lem:operL} 
\be
|v|_{\om,3} = O \Big( |v|_\om + |f(v+ w(v))|_{\om} \Big) = 
O \Big( |v|_\om + |v|_\om^p \Big( 
1 + \frac{|v|_\om^{p-1}}{\gamma} \Big) \Big) =
O(|v|_\om )
\ee
because $ v \in {\cal D}_\rho $. Hence, by remark \ref{rem:red},
$|v|_{\om,3}^{p-1} \slash \gamma \leq 1 \slash 2 \wtilde{C} $
provided the constants $C_{10}$, $C_{11}$, $C_{12}$, $C_{13}$ 
have been chosen small enough.

\section{Appendix}

\begin{pfn}{\sc of lemma} \ref{composition}.
By standard results, if $ u\in L^{\infty} \cap H^1 ( {\bf T} \times
(0,\pi))$, then $f(u) \in H^1 ({\bf T} \times
(0,\pi))$ and  $(f(u))_x \equiv f'(u) u_x $, 
$(f(u))_t \equiv f'(u) u_t $. We get  
\be\label{eq:Sobol}
| (f(u))_x |_{L^2} = \Big( \int_\Omega |(f(u))_x|^2 \Big)^{1/2} = 
\Big(  \int_\Omega |f'(u)|^2 |u_x|^2 \Big)^{1/2} 
\leq |f'(u)|_\infty ||u||_{H^1} \leq C |u|_\infty^{p-1} ||u||_{H^1}
\ee
for $ | u |_\infty $ small enough.
An analogous estimate holds for $ | ( f(u))_t |_{L^2} $. 
Since $| f (u) | \leq C |u |_\infty^p$, 
we obtain (\ref{nemitski}).

We now prove that $u \to f(u) $ is $C^1(X, X)$. It is easy 
to show that it is Gateaux differentiable and that
$D_G f(u)[h] = f'(u)h $. Moreover with estimates similar as before
we have that $u \to D_G f(u)$ is continuous. Hence the Nemitski operator
$u \to f(u)$  is in $C^1 ( X, X ) $ and its Frechet differential is
$h \to Df(u)[h]= D_G f(u)[h] = f'(u) h $.
Estimate (\ref{nemitski2}) 
can be obtained in the same way as  (\ref{nemitski}).
\end{pfn}

\begin{pfn}{\sc of lemma} \ref{lem:operL}. 
Writing $ w(t,x) = \sum_{l \geq 0, j \geq 1, j \neq l} 
w_{l,j} \cos (lt) \sin (jx) \in X $ 
we have that 
$$
L_\om^{-1} w (t,x) = \sum_{l\geq 0, j \geq 1, j \neq l}
\frac{w_{l,j}}{(\om l-j)(\om l+j)} \cos (lt) \sin (jx). 
$$
Since $\om \in W_\gamma$, it is clear that 
$||L_\om^{-1}w||= O(||w||/\gamma)$.

We now prove that for all $ w \in W $
\be\label{eq:S}
S := \sum_{ l \geq 0, j \geq 1, j \neq l} 
\frac{|w_{l,j}|}{|\om l-j|(\om l+j)} \leq C \Big( 
|w|_{L^2} + \frac{|\om-1|^{1/2}}{\gamma} ||w||_{H^1} \Big).
\ee
Note that (\ref{eq:S}) implies 
\be\label{L-1w}
|L_\om^{-1} w |_{\infty} \leq S \leq  C \Big( 
|w|_{L^2} + \frac{|\om-1|^{1/2}}{\gamma} ||w||_{H^1} \Big).
\ee
For $l \in {\bf N}$, let $ e(l) \in {\bf N} $ be defined by
$ |e(l)- \om l|= \min_{j \in {\bf N}} |j-\om l|.$
Since $ \om $ is not rational, 
$e(l)$ is the only integer $e$ such that
$|e-\om l|<1/2$.  Let us define
$$
S_1=\sum_{l\geq 0, j \geq 1, j \neq l, j \neq e(l)} 
\frac{|w_{l,j}|}{|\om l-j|(\om l+j)}
$$
and 
$$
S_2=\sum_{ l \geq 0, e(l) \neq l} 
\frac{|w_{l,e(l)}|}{|\om l-e(l)|(\om l+e(l))}.
$$
For $j \neq e(l)$ we have that
\be\label{est:smde}
|\om l - j | \geq  \frac{|j-e(l)|}{2} \qquad {\rm and } 
\qquad \om l + j \geq \frac{ | j - e(l) | + l }{4}.
\ee
Indeed, for $j \neq e(l)$ we have that 
$\ |j -\om l | \geq |j-e(l)|-|e(l)-\om l| \geq  |j-e(l)|-1/2
\geq |j-e(l)|/2$.
Moreover, since $ | e(l) - \om l | < 1 \slash 2 $, 
it is easy to see (remember that $\om \geq 1/2$) that 
$ e(l) + l \leq 4 \om l $ and hence  
$\ |j-e(l)|+l \leq j+e(l)+l  \leq 4
(j+\om l)$. Defining $w_{l,j}$ by $w_{l,j}=0$ 
if $j \leq 0 $ or $j=l$, we get from (\ref{est:smde}) 
$$
S_1 \leq \sum_{ l \geq 0, j \in {\bf Z}, j \neq e(l)} 
\frac{8|w_{l,j}|}{|j-e(l)|(|j- e(l)|+l)}.
$$ 
Hence, by the Cauchy-Schwarz inequality, $S_1 \leq 8R_1
|w|_{L^2}$, where
$$
R_1^2= \sum_{l\geq 0, j\in {\bf Z}, j\neq e(l)} 
\frac{1}{(j-e(l))^2(|j- e(l)|+l)^2} = 
\sum_{l\geq 0, j\in {\bf Z}, j\neq 0} 
\frac{1}{j^2(|j|+l)^2} \leq \sum_{l\geq 0, j\in {\bf Z}, j\neq 0} 
\frac{1}{j^2(1+l)^2}  < \infty.
$$
In order to find an upper bound for $S_2$, we observe that
if $ l < $ $ 1 / (2 |\om-1|) $ then $ |\om l-l| < 1/ 2 $ and then $ e(l)=l$.
Hence  we can write
$$
S_2=\sum_{ l\geq 1/2|\om-1|, e(l) \neq l} 
\frac{|w_{l,e(l)}|}{|\om l-e(l)|(\om l+e(l))}.
$$ 
By the properties of $ \om \in W_\gamma $, for $ l \neq e(l)$,
$ |\om l - e(l)||\om l + e(l) | \geq $  $ \gamma (\om l + e(l)) 
l^{-1} \geq$ $ \gamma  $. Hence, still by the
Cauchy-Schwarz inequality,
$$S_2 \leq \frac{1}{\gamma} \sum_{ l\geq 1/2|\om-1|, e(l) \neq l}
|w_{l,e(l)}|  \leq \frac{ 2 R_2 }{ \gamma } || w ||_{H^1},$$ 
where
$$
R_2^2=\sum_{ l\geq 1/2|\om-1|} \frac{1}{l^2} \leq 2 |\om-1|.
$$
The estimates $ | L^{-1}_\om w |_\om \leq (C \slash \gamma ) | w |_\om $ 
$ \forall w \in W $ and 
$| L^{-1}_\om \Pi_W  u|_\om \leq 
(C_3 \slash \gamma ) |u|_\om $ $ \forall u \in X $ 
are an immediate consequence of the bound of 
$||L_\om^{-1}w||$ and (\ref{L-1w}), since
$|\Pi_W u|_{L_2} \leq |u|_{L_2}$ and
$||\Pi_W u|| \leq ||u||$
For the estimate of $|L^{-1}_\om w|_{L^2}$ ($w\in W$), we use
the decomposition $|L^{-1}_\om w|_{L^2}^2=S'_1+S'_2$, where
$$
S'_1 := \sum_{j>0, l\geq 0, j\neq l, j\neq e(l)} 
\frac{w_{l,j}^2}{(\om l-j)^2(\om l+j)^2},  \quad 
\quad S'_2 := \sum_{ l\geq 1/2|\om-1|, e(l) \neq l} 
\frac{w_{l,e(l)}^2}{(\om l-e(l))^2(\om l+e(l))^2}.
$$
Recall that if
$ l < 1/2 |\om-1| $ then $e(l)=l$.
Arguing as before it is obvious that $ S'_1 \leq C |w|_{L^2}^2$. We have
$$
S'_2 \leq \sum_{ l\geq 1/2|\om-1|, e(l) \neq l}
\frac{w_{l,e(l)}^2}{\gamma^2} \leq \frac{4(\om-1)^2}{\gamma^2}
\sum_{ l\geq 1/2|\om-1|, e(l) \neq l} 
l^2 w_{l,e(l)}^2 = O \Big(\frac{\e^2}{\gamma^2}  ||w||_{H^1}^2 
\Big) = O \Big( \frac{\e}{\gamma^2}  |w|_\om^2 \Big),
$$
which yields  estimate (\ref{estsupp}).
We now prove (\ref{diffLLom}). 
Writing $ r(t,x) = \sum_{ l \geq 0, j \geq 1} 
r_{l,j} \cos (lt) \sin (j x) $, $ s(t,x) = 
\sum_{l \geq 0, j \geq1 } 
s_{l,j} \cos (lt) \sin (j x) $ we get 
$$
L_1^{-1} \Pi_W s= 
\sum_{j \neq l} \frac{s_{l,j}}{l^2 - j^2} \cos (lt)
 \sin (j x), \qquad
 \qquad 
L_\om^{-1} s = \sum_{j \neq l} 
\frac{s_{l,j}}{\om^2 l^2 - j^2}  \cos (lt) \sin (j x)
$$
and
$$
S:=\int r (L_\om^{-1} - L^{-1})(\Pi_W s) \ dt \ dx = 
\pi^2 \sum_{j\neq l} 
\frac{s_{l,j}r_{l,j}\ \e l^2}{(\om^2l^2-j^2)(l^2-j^2)}=\pi^2 (S_1+S_2),
$$
where
$$
S_1 :=\sum_{j\neq l, j\neq e(l)} 
\frac{s_{l,j}r_{l,j}\ \e l^2}{(\om l-j)(\om l +j)(l-j)(l+j)},
\qquad 
S_2 :=\sum_{l\neq e(l)} 
\frac{s_{l,e(l)}r_{l,e(l)}\ \e l^2}{(\om l-e(l))(\om l +e(l))
(l-e(l))(l+e(l))}.
$$
For $j\neq e(l)$ and $j\neq l$ we have that 
$| \om l - j | \geq |j - e(l)| \slash 2 \geq 1 \slash 2 $ hence
$| (\om l-j)(\om l
+j)(l-j)(l+j) | \geq (1/2) \om l^2$.  We obtain
$$
|S_1| \leq \frac{2\e}{\om} \sum_{j\neq l}|r_{l,j}| |s_{l,j}| =O
\Big( \frac{2\e}{\om} | r |_{L^2} |s |_{L^2} \Big)
=O\Big( \frac{2\e}{\om} |r|_\om 
|s|_\om \Big). 
$$ 
As before, for $\om \in W_\gamma $ and $ l \neq e(l) $, 
we have $|\om l - e(l)| \geq \gamma \slash l $ and so 
$| (\om l-e(l))(\om l +e(l))
(l-e(l))(l+e(l) | \geq  \gamma l$. Moreover $ l\neq e(l)$ 
implies $ l \geq 1/(2|\om -1|)$. Hence
$$
|S_2| \leq  C \frac{\e^2}{ \gamma} \sum_{l \geq 1/(2\e)}
|r_{l,e(l)}|  |s_{l,e(l)}| l^2 =O\Big( \frac{\e^2}{ \gamma}
|| r ||_{H^1} || s ||_{H^1} \Big).   
$$
Estimate (\ref{diffLLom}) follows straightforwardly from the
above estimates of $S_1$ and $S_2$.
Finally (\ref{Lomest}) is an immediate consequence of (\ref{diffLLom})
since, for all $(r,s) \in X \times X$,
$$
\int_{\Om} r L^{-1}(\Pi_W s) = O( |r|_{L^2} |s|_{L^2})= O(|r|_\om |s|_\om).
$$
\end{pfn}

\begin{pfn}{\sc of lemma} \ref{lem:red}.
It is sufficient to show that the operator ${\cal G} : W \to W $ defined by 
${\cal G} (w) :=  L_\om^{-1} \Pi_W f(v+w)$ is a contraction on a ball
$ B_\d \subset W $ for some $ \d $ if
$|v|_\om^{p-1}/\gamma$ is small enough. 
Since, by lemma \ref{lem:operL},
$|| L_\om^{-1} \Pi_W ||_\om \leq C_3 / \gamma $ and by 
(\ref{nemitski}), we have that, for all $w\in B_\d$ 
\be\label{ballball}
|{\cal G} (w) |_\om \leq \frac{C_3}{\gamma} |f(v+w)|_\om  \leq
\frac{C}{\gamma} |v+w|_\om^p \leq \frac{C}{\gamma} 
\Big( |v|_\om + \delta \Big)^p 
\leq \frac{\ov{C}}{\gamma} \Big( | v|_\om^p  + \delta^p \Big).
\ee
for some $ \ov{C} > 0 $.
Moreover by lemma \ref{composition} the  operator ${\cal G}$  is
differentiable and $D{\cal G}(w)[h] = L_\om^{-1} \Pi_W
(f'(v+w)h)$. Hence by (\ref{L-1w}),
\be \label{DGest}
|D{\cal G}(w)[h]|_\om \leq \frac{C_3}{\gamma} |f'(v+w)h|_\om \leq
 \frac{C}{\gamma} |v+w|_\om^{p-1} |h|_\om \leq
\frac{\wtilde{C}}{\gamma}  \Big( |v|_\om^{p-1}+ \delta^{p-1} \Big) |h|_\om
\ee   
for some $\wtilde{C} > 0$.
Let us choose $ \delta := |v|_\om$ and
$ \rho := (4\max(\ov{C},\wtilde{C}) )^{-1}$. 
For $v \in {\cal D}_\rho$, 
$(\ov{C} / \gamma) (|v|_\om^p +\delta^p) \leq \delta /2$ and
$(\wtilde{C}/ \gamma)  (|v|_\om^{p-1}+ \delta^{p-1})\leq 1/2$. 
Hence ${\cal G}$ maps $B_\delta$ into the closed ball $\ov{B}_{\delta/2}$  
and is a contraction on ${B}_{\delta}$.
By the Contraction Mapping Theorem there is a unique $w(v)$
in ${B}_{\delta}$ such that ${\cal G}(w(v))=w (v)$. 
We have 
\be\label{eq:stimco}
|w(v)|_\om =| {\cal G} (w(v))|_\om  \leq \frac{2\ov{C}}{\gamma}|v|_\om^p
\ee
by (\ref{ballball}) and the choice of $ \delta $.
Moreover
$$
| w(v) |_{L^2}=| L_\om^{-1} \Pi_W f( v + w ( v ) ) |_{L^2} \leq
C \Big( 1 + \frac{\sqrt{\e}}{\gamma} \Big) | f ( v + w ) |_\om
\leq C \Big( 1+ \frac{\sqrt{\e}}{\gamma} \Big) |v|_\om^p,
$$
by (\ref{estsupp}) and (\ref{eq:stimco}).
In addition
\begin{eqnarray*}
\Big| \int_{\Om} r \Big( w(v) -L^{-1}_\om \Pi_W f(v) \Big) \Big| & = & 
\Big| \int_\Om r \ L_\om^{-1}\Pi_W (f(v+w(v))-f(v)) \Big| \\
&\leq& C_4 \Big( 1+\frac{\e}{\gamma} \Big) |f(v+w(v))-f(v)|_\om |r|_\om \\
&\leq& C \Big( 1+\frac{\e}{\gamma} \Big) ( |v|_\om^{p-1} |w(v)|_\om) |r|_\om \\
&\leq& \frac{C}{\gamma} \Big( 1+\frac{\e}{\gamma} \Big) 
|v|_\om^{2p-1} |r|_\om, 
\end{eqnarray*}
by (\ref{Lomest}) and the bound on $| w ( v ) |_\om $ 
given in (\ref{eq:stimco}). This proves $ii)$.

In order to prove $iii$), let us call ${\cal I} : X\to X$ the linear
operator defined by ${\cal I}u (t,x)=u(t+\pi, \pi-x)$.
It is easy to see that ${\cal I}$ and $L_\om^{-1}\Pi_W$ commute.
Hence, since   $ - v = {\cal I}v$, 
$$
w(v) =  L_\om^{-1} \Pi_W f(v + w(v) ) \ \Rightarrow
\ {\cal I} w(v)=  
L_\om^{-1} \Pi_W {\cal I}f( v + w(v)) =
L_\om^{-1} \Pi_W f( -v + {\cal I}w(v)).
$$
This implies, by the uniqueness of the solution of $w=
L_\om^{-1} \Pi_W f( -v + w)$ in the ball of radius $|v|_\om$,
that $ w(v)( t + \pi, \pi -x ) = w(-v)( t, x) $.

For $iv)$, we remark that if $v\in V_n$, then 
$w(v)$ is the unique solution in an appropriate ball 
of the equation ${\cal G}(w)=w$, where the map ${\cal G}$  
satisfies ${\cal G}(W_n) \subset W_n$ ($W_n$ is the closed subspace
of $W$ containing 
 $w\in W$ which are $2\pi /n$ periodic in time). 
Since ${\cal G}^k (0) \to w(v) $ in $W$, this implies that
$w(v) \in W_n$. 

Finally, the map $ (v,w) \to {\cal G}(v; w):= w-L_\om^{-1}
\Pi_W f(v+w)$ is of class $C^1 $  and its differential with 
respect to $w$ is invertible at any point $(v,w(v))$ by
 lemma \ref{composition} and the previous bounds. 
As a consequence of the Implicit Function 
Theorem the map $ v \to w(v)$ is in $C^1({\cal D}_\rho,W)$.
\end{pfn}

\begin{pfn} {\sc of lemma} \ref{lemexpl}.
We define $ \alpha $, $ a_1 $, $ a_2 $, $ \wtilde{m} $ by 
$ \alpha= \la m\ra $,
$$
a_1(s_1) := \la m(s_1, \cdot )-\alpha\ra_{s_2}= \frac{1}{2\pi}
\int_0^{2\pi} m(s_1,s_2)-\alpha \ ds_2, \ \
a_2(s_2) := \la m( \cdot ,s_2)-\alpha\ra_{s_1}= \frac{1}{2\pi}
\int_0^{2\pi} m(s_1,s_2)-\alpha \ ds_1
$$
and $\wtilde{m} (s_1,s_2) := m(s_1,s_2) - a_1(s_1)- a_2(s_2) - \alpha $. 
It is straightforward to check that $\la a_i \ra =0$,
$\la \wtilde{m} \ra_{s_2} (s_1)=0$, $ \la \wtilde{m} \ra_{s_1}
(s_2)=0$.  We have to prove that $ a_1=a_2 $. Since $ w\in W$,      
for all $p$ odd and $2\pi$ periodic,
$$
\int w(t,x) (p(t+x)-p(t-x)) \ dt \ dx = \frac{1}{2} 
\int_{{\bf T}^2} m(s_1,s_2) (p(s_1)-p(s_2)) \ ds_1 \ ds_2 =0. 
$$
By the definition of $a_i$, we obtain that for all $p$ odd and
$2\pi$-periodic
$$
\int_{\bf T} (a_1(s)-a_2(s) ) p(s) \ ds =0.
$$
Therefore $a_1-a_2$ is even. Now $ w(t,x) = w(-t,x) $ for all $ t\in
{\bf T}$ and for all $ x \in [0,\pi] $, and this implies that 
$m(-s_2,-s_1)=m(s_1,s_2)$. As a consequence $ a_2 ( s ) = a_1 ( - s ) $, so
that $ a_1-a_2 $ is odd. Since $a_1-a_2$ is both odd and even, 
$a_1=a_2$.

We now turn to the expression of $\int w L^{-1} (w)$. Let us define
$\beta_1$, $\beta_2$, $\beta_3$ on $[0,2\pi]^2$ by
$$
\begin{array}{rcl}
\beta_1(s_1,s_2)&:=& - M(s_1,s_2) +
\frac{1}{2}(M(s_1,s_1)+M(s_2,s_2)), \\
\beta_2(s_1,s_2)&:=& (s_1-s_2) (A(s_1)-A(s_2)), \\
\beta_3(s_1,s_2)&:=&- \dps \frac{\alpha}{8} (s_1-s_2)(2\pi -(s_1-s_2)),
\end{array}
$$ 
and $B_1,B_2,B_3$ on $[0,2\pi] \times [0,\pi]$ by
$B_i(t,x) := \beta_i(t+x,t-x)$. As an immediate consequence of its
definition, $B_i$ is $2\pi$-periodic w.r.t. $t$ and for $x\in \{0,
\pi\}$, $B_i(t,x)=0$. Hence $B_i \in X$. Moreover
$-\partial_{s_1} \partial_{s_2} \beta_i = (1/4) \lambda_i$, where
$$
\lambda_1(s_1,s_2):=\wtilde{m}(s_1,s_2), \qquad \lambda_2(s_1,s_2) :=
a(s_1)+a(s_2), \qquad \lambda_3(s_1,s_2) := \alpha.
$$
This implies that $[-(\partial_t)^2 + (\partial_x)^2] B_i=L_i$ where we
have defined $ L_i(t,x) := \lambda_i(t+x, t-x)$. As a result
$ B_1+B_2+B_3 \in L^{-1} w + V$ and
$$
\int_{[0,2\pi] \times [0,\pi]} w(t,x) L^{-1} (w) (t,x) \ dt \ dx = 
\int_{[0,2\pi] \times [0,\pi]}
(B_1+B_2+B_3)(t,x) (L_1+L_2+L_3)(t,x) \ dt \ dx.  
$$
We observe that, since $B_i$ vanishes for $x=0$ and $x=\pi$,  for all
$i,j$, 
$$
\int B_i  L_j \ dt \ dx =
\int B_i \Big( (B_j)_{tt}-(B_j)_{xx} \Big) \
dt \ dx =\int
\Big( (B_i)_{tt}-(B_i)_{xx} \Big) B_j \ dt \ dx = \int
 L_i  B_j \ dt \ dx.   
$$ 
Hence 
$$
\int w L^{-1} w= \int_{[0,2\pi] \times [0,\pi]} B_1 L_1
+2\int  B_1 L_2 + 2\int  B_1
L_3 + \int  B_2 L_2 + 2 \int  B_2
L_3 + \int  B_3 L_3.
$$
We have 
$$
\int  B_1 L_1=\frac{1}{2} \int_{{\bf T}^2}
(-M(s_1,s_2)+ \frac{1}{2} M(s_1,s_1) + \frac{1}{2} M(s_2,s_2))
\wtilde{m} (s_1,s_2) \ ds_1 \ ds_2 =- \frac{1}{2} \int_{{\bf T}^2}
M(s_1,s_2) \wtilde{m}(s_1,s_2) \ ds_1 \ ds_2 
$$
because $\la \wtilde{m} \ra_{s_i}=0$. 
$$
\int B_1 L_2 = \frac{1}{2} \int_{{\bf T}^2}
(-M(s_1,s_2)+ \frac{1}{2} M(s_1,s_1) + \frac{1}{2} M(s_2,s_2))
(a(s_1)+a(s_2)) \ ds_1 \ ds_2 =\pi \int_{{\bf T}} M(s,s) a(s) \ ds
$$
because $\la M \ra_{s_i}=0$ and $\la a \ra=0$.
$$
\int B_1 L_3=\frac{\alpha}{2} \int_{{\bf T}^2}
\Big( -M(s_1,s_2)+ \frac{1}{2} M(s_1,s_1) + \frac{1}{2} M(s_2,s_2) \Big)=-\pi
\alpha \int_{{\bf T}} M(s,s) \ ds 
$$ 
because $\la M \ra=0$. 
\begin{eqnarray*}
\int  B_2 L_2 &=& \int_{0}^{2\pi}
\int_0^{\pi} 2x \Big( A(t+x)-A(t-x) \Big) \Big( 
a(t+x)+a(t-x) \Big) \ dt \ dx \\
&=& \int_0^{\pi}  dx  \ 8x \int_0^{2\pi} \Big( A(t+x)-A(t-x) \Big)
\Big( A'(t+x)+A'(t-x) \Big) \
dt \\ 
&=& \int_0^{\pi}  dx  \ 8x \int_0^{2\pi} \frac{d}{dt}
\Big[\frac{(A(t+x)-A(t-x))^2}{2} - A(t-x)^2\Big] + 2A(t+x) A'(t-x) \ dt
\\
&=& \int_0^{\pi}  dx  \ 16x \int_0^{2\pi} A(t+x)A'(t-x) \ dt =
\int_0^{\pi}  dx  \ 16x \int_0^{2\pi} A(s) A'(s-2x) \ ds \\
&=& \int_0^{2\pi} ds \  A(s) \int_0^{\pi} 16x A'(s-2x) \ dx =
\int_0^{2\pi} ds \  A(s) \Big[ 8x A(s-2x) \Big]_{x=\pi}^{x=0} \\
&=& -8\pi \int_0^{2\pi} A(s)^2 \ ds. 
\end{eqnarray*}
In the fore-last line, we have integrated by parts (w.r.t. $x$) and used
the fact that $\la A \ra =0$.
$$
\int  B_2 L_3 = \alpha \int_0^{\pi} dx \
2x \int_0^{2\pi} (A(t+x)-A(t-x)) \ dt =0,
$$
still because $\la A \ra=0$. At last
$$
\int B_3 L_3 =-\frac{\alpha}{8}
\int_0^{2\pi} \int_0^{\pi} 2x (2\pi -2x) \ dt \ dx =
-\frac{\alpha^2 \pi^4}{6}
$$
Summing up, we get (\ref{calculpenible}).
\end{pfn}


\begin{thebibliography}{10}

 \bibitem{AB}
 A. Ambrosetti, M. Badiale,
 {\it Homoclinics: Poincar\'e-Melnikov type results
 via a variational approach},
 Annales I. H. P. - Analyse nonlin., vol. 15, n.2, pp. 233-252, 1998.

\bibitem{ACE} A. Ambrosetti, V. Coti-Zelati, I. Ekeland,
{\it Symmetry breaking in Hamiltonian systems},
Journ. Diff. Eq. 67, 165-184, 1987.

\bibitem{AM} A. Ambrosetti, A. Malchiodi, {\it A multiplicity result for 
the Yamabe problem on $ S^n $}, J. Funct. Anal. 168, 529-561, 1999.

\bibitem{AP} A. Ambrosetti, G. Prodi, {\it A primer of nonlinear
analysis}, Cambridge University Press, Cambridge, 1993.

\bibitem{AR} A. Ambrosetti, P. Rabinowitz:
{\it Dual Variational Methods in Critical Point Theory and Applications}, 
Journ. Func. Anal, 14, 349-381, 1973.

\bibitem{B1} D. Bambusi: {\it Lyapunov Center Theorems for some nonlinear
PDE's: a simple proof}, Ann. Sc. Norm. Sup. di Pisa, Ser. IV, 
vol. XXIX, fasc. 4, 2000.

\bibitem{BP1} D. Bambusi, S. Paleari, {\it Families of periodic solutions 
of resonant PDE's}, J. Nonlinear Sci., 11, pp. 69-87, 2001.

\bibitem{BP2} D. Bambusi, S. Paleari, {\it Families of periodic orbits
for some PDE's in higher dimensions}, Comm. Pure and Appl. Analysis,
Vol. 1, n.4, 2002.

\bibitem{Bar} T. Bartsch, 
{\it A generalization of the Weinstein-Moser Theorems on periodic 
orbits of a Hamiltonian system near an equilibrium},
Ann. Inst. Henri Poincar\'e, Anal. Non Lin\'eaire 14, No.6, 691-718, 1997.

\bibitem{BM} M. Berti, A. Malchiodi, 
{\it Multiplicity results and multibump solutions for the 
Yamabe problem on $S^n$"}, J. Funct. Anal., 
180, 1, 2001. 

\bibitem{BL} G.D. Birkhoof, D.C. Lewis, {\it On the periodic motions
near a given periodic motion of a dynamical system}, Ann. Mat. 
12, 117-133, 1933.

\bibitem{B0} J. Bourgain, {\it Construction of Quasi-Periodic solutions
for Hamiltonian perturbations of linear equations 
and applications to Nonlinear PDE}, Int. Math. Res. Not, n.11, 1994. 

\bibitem{Bo1} J. Bourgain, {\it Construction of periodic solutions 
of nonlinear wave equations in higher dimension}, Geom. and Func. Anal.,
vol. 5, n.4, 1995.

\bibitem{B3} J. Bourgain, {\it Quasi-periodic solutions of Hamiltonian 
perturbations of $2D$ linear Schrodinger equations}, 
Ann. of Math., 148,  363-439, 1998.

\bibitem{BCN} H. Brezis, J. M. Coron, L. Nirenberg, {\it Free
vibrations for a non-linear wave equation and a Theorem of 
P. Rabinowitz}, Comm. Pure and Appl. Math. 31, 1-30, 1978.

\bibitem{Co} J. M. Coron, {\it Periodic solutions of a Non-linear 
wave equation without assumption of monotonicity}, Math. Ann. 
262, 273-285, 1985.

\bibitem{CW} W. Craig, E.Wayne, {\it Newton's method and periodic solutions 
of nonlinear wave equation}, Comm. Pure and Appl. Math, 
vol. XLVI, 1409-1498, 1993.

\bibitem{DST} L. DeSimon, G. Torelli, {\it Soluzioni periodiche 
di equazioni alle derivate parziali di tipo iperbolico nonlineari},
Rend. Sem. Mat. Univ. Padova, 40,  380-401, 1968.


\bibitem{FR} E. R. Fadell, P. Rabinowitz, {\it Generalized 
cohomological index theories for the group actions 
with an application to bifurcations question for Hamiltonian
systems}, Inv. Math. 45, 139-174, 1978.


\bibitem{K2} S. B. Kuksin, 
{\it Perturbation of conditionally periodic solutions 
of infinite-dimensional Hamiltonian systems},
Izv. Akad. Nauk SSSR, Ser. Mat. 52, No.1, 41-63, 1988.

\bibitem{Ku} S. B. Kuksin, {\it Analysis of Hamiltonian PDE's}
Oxford Lecture Series in Mathematics and its Applications. 19. 
Oxford University Press, 2000.

\bibitem{Ly} A. M. Lyapunov, {\it Probl\`eme g\'en\'eral de la  stabilit\'e
du mouvement}, Ann. Sc. Fac. Toulouse, 2, 203-474, 1907.

\bibitem{LS} B. V. Lidskij, E. I. Shulman, {\it Periodic
solutions of the equation $u_{tt} - u_{xx} + u^3  = 0 $},
Funct. Anal. Appl., 22, 332-333, 1988. 

\bibitem{Mo1} J. Moser, {\it Proof of a generalized form of a fixed 
point theorem due to G. D. Birkhoof}, Geometry and Topology,
Lectures Notes in Math. 597, 464-494, 1976.

\bibitem{Mo} J. Moser, {\it Periodic orbits near an Equilibrium 
and a Theorem by Alan Weinstein}, Comm. on Pure and Appl. Math.,
vol. XXIX, 1976.

\bibitem{Po} H. Poincar\'e, {\it Les m\'ethodes nouvelles de la m\'ecanique
c\'eleste}, Paris, 1899.

\bibitem{R1} P. Rabinowitz, {\it Time periodic solutions
of nonlinear wave equations},  
Manusc. Math. 5, 165-194, 1971.

\bibitem{R0} P. Rabinowitz, {\it Periodic solutions of nonlinear 
hyperbolic partial differential equations}, Comm. Pure Appl. Math.,
20,  145-205, 1967.

\bibitem{R} P. Rabinowitz, {\it Free vibration of a seminilear 
wave equation}, Comm. on Pure Appl. Math, 31, 
31-68, 1978.


\bibitem{Str} M. Struwe, {\it  Variational methods. 
Applications to nonlinear Partial Differential Equations and
Hamiltonian Systems}, Springer Verlag, 1990.


\bibitem{W1} E. Wayne, {\it Periodic and quasi-periodic solutions 
of nonlinear wave equations via KAM theory},
Commun. Math. Phys. 127, No.3, 479-528, 1990.

\bibitem{We} A. Weinstein, {\it Normal modes for 
Nonlinear Hamiltonian Systems}, Inv. Math, 20, 47-57, 1973. 

\end{thebibliography}
\end{document}